\documentclass[10pt]{amsart}
\usepackage[pagebackref]{hyperref}

\renewcommand{\d}{\mathrm{d}}

\newcommand{\ts}{\textstyle }
\newcommand{\ds}{\displaystyle }

\newcommand{\cE}{{\mathcal E}}

\newcommand{\cX}{{\mathcal X}}

\newcommand{\bbR}{{\mathbb R}}
\newcommand{\R}[1]{\ensuremath{\mathbb R^{\,#1}}} 

\newcommand{\E}[1]{\ensuremath{\mathbb E^{\,#1}}} 

\newcommand{\bbZ}{{\mathbb Z}}

\newcommand{\ot}{\otimes}

\newcommand{\fff}{\mathrm{I}}
\newcommand{\sff}{\mathrm{I\!I}}

\newcommand{\SO}{\operatorname{SO}}

\newcommand{\Or}{\operatorname{O}}

\newcommand{\End}{\operatorname{End}}
\newcommand{\id}{\operatorname{id}}

\newcommand{\ab}{\mathbf{a}}
\newcommand{\bb}{\mathbf{b}}
\newcommand{\eb}{\mathbf{e}}
\newcommand{\nb}{\mathbf{n}}
\newcommand{\tb}{\mathbf{t}}
\newcommand{\ub}{\mathbf{u}}

\newcommand{\xb}{\mathbf{x}}
\newcommand{\yb}{\mathbf{y}}
\newcommand{\wb}{\mathbf{w}}
\newcommand{\zb}{\mathbf{z}}

\newcommand{\w}{{\mathchoice{\,{\scriptstyle\wedge}\,}{{\scriptstyle\wedge}}
      {{\scriptscriptstyle\wedge}}{{\scriptscriptstyle\wedge}}}}

\numberwithin{equation}{subsection}

\newtheorem{theorem}{Theorem}
\newtheorem{lemma}{Lemma}
\newtheorem{proposition}{Proposition}
\newtheorem{corollary}{Corollary}

\theoremstyle{remark}
\newtheorem{definition}{Definition}

\newtheorem{remark}{Remark}
\newtheorem{example}{Example}

\begin{document}

\author[R. Bryant]{Robert L. Bryant}
\address{Duke University Mathematics Department\\
         P.O. Box 90320\\
         Durham, NC 27708-0320}
\email{\href{mailto:bryant@math.duke.edu}{bryant@math.duke.edu}}
\urladdr{\href{http://www.math.duke.edu/~bryant}%
         {http://www.math.duke.edu/\lower3pt\hbox{\symbol{'176}}bryant}}

\dedicatory{This article is dedicated to Shiing-Shen Chern, 
whose beautiful works and gentle encouragement have had 
the most profound influence on my own research.}

\title[Prescribed shape operator]
      {On surfaces with  \\
          prescribed shape operator}

\date{July 20, 2001}

\begin{abstract}
The problem of immersing a simply connected surface with 
a prescribed shape operator is discussed.  From classical 
and more recent work (see~\cite{Fe:survey} for a survey), 
it is known that, 
aside from some special degenerate cases, such as when 
the shape operator can be realized by a surface with one family 
of principal curves being geodesic, the space of such
realizations is a convex set in an affine space of dimension
at most~$3$.  The cases where this maximum dimension of 
realizability is achieved have been classified and it is known that
there are two such families of shape operators, one
depending essentially on three arbitrary functions of one variable
(called Type~I in this article) and another depending essentially 
on two arbitrary functions of one variable (called Type~II in this article).

In this article, these classification results are rederived,
with an emphasis on explicit computability of the space of solutions.
It is shown that, for operators of either type, their realizations
by immersions can be computed by quadrature.  Moreover, explicit
normal forms for each can be computed by quadrature together with, 
in the case of Type~I, by solving a single linear second order ODE 
in one variable. (Even this last step can be avoided in most 
Type~I cases.) 

The space of realizations is discussed in each case, 
along with some of their remarkable geometric properties.
Several explicit examples are constructed (mostly already
in the literature) and used to illustrate various features of the problem.
\end{abstract}



\thanks{
The research for this article was begun during a June 2001 
conference at the Mathematisches Forschungsinstitut at Oberwolfach.
The author thanks the MFO for its hospitality.
Thanks also to Duke University for its support 
via a research grant 
and to the National Science Foundation 
for its support via DMS-9870164.\hfill\break
\hspace*{\parindent} 
This is Version~$3.0$. The most recent version 
can be found at arXiv:math.DG/0107083 .
}

\maketitle

\setcounter{tocdepth}{2}
\tableofcontents

\section{Introduction}\label{sec: intro}

\subsection{The fundamental forms}\label{ssec: fund forms}
In classical surface theory in Euclidean space, given
an immersion~$\xb:D\to\E3$ of a surface~$D$ into Euclidean $3$-space,
one can define its first fundamental form~$\fff_\xb=\d\xb\cdot\d\xb>0$ 
and, given a choice of unit normal~$\nb:D\to S^2$ for~$\xb$ 
(i.e., $\nb\cdot\d\xb=0$ and~$\nb\cdot\nb=1$), one can define its
second fundamental form~$\sff_\xb = -\d\nb\cdot\d\xb$.  The
quantities~$\fff_\xb$ and~$\sff_\xb$ are unchanged if one 
composes~$\xb$ with an isometry of~$\E3$ and replaces~$\nb$ by
its corresponding image under this isometry.  Moreover, two normally
oriented immersions~$\xb,\yb:D\to\E3$ agree up to isometry 
to second order at a point~$p\in D$ if and only 
if~$\fff_\xb(p)=\fff_\yb(p)$ and~$\sff_\xb(p)=\sff_\yb(p)$.  
Thus, the two quadratic forms~$(\fff_\xb,\sff_\xb)$ contain all 
of the second-order information about a normally oriented 
immersion~$\xb$ that is invariant under Euclidean isometries.

\subsection{Bonnet's theorem and rigidity}\label{ssec: bonnet and rigidity}
One of the most classical theorems in the subject is Bonnet's theorem,
which asserts that a given pair of quadratic forms~$(\fff,\sff)$ defined
on a simply connected surface~$D$ can be realized by an immersion~$\xb$
with choice of unit normal~$\nb$ if and only if~$\fff$ is positive
definite and the pair~$(\fff,\sff)$ 
satisfy the Gau{\ss} and Codazzi equations.
Moreover, $\xb$ and~$\nb$ (when they exist) are unique up to an 
isometry of~$\E3$.   Since specifying a pair of quadratic forms on
a surface is tantamount to choosing six arbitrary functions of two 
variables while choosing an immersion of the surface into~$\E3$ is
tantamount to choosing three arbitrary functions of two variables,
it is not surprizing that there exist such compatibility conditions
on pairs~$(\fff,\sff)$ in order that they be realizable by an
immersion.

\subsection{Isometric embedding}\label{ssec: iso embed}
It is natural to look at problems that are not as overdetermined as
Bonnet's.  For example, the problem of finding an immersion~$\xb$ that
realizes a given positive definite quadratic form~$\fff>0$ as its
first fundamental form is known as the \emph{isometric embedding
problem} and has a long history in differential geometry.  
The equation~$\fff_\xb=\fff$, regarded as an equation for~$\xb$,
is determined in the na\"{\i}ve sense (i.e., it is three equations
for three unknowns), but its behavior is rather subtle.  It is known
to be locally solvable when~$\fff$ is real-analytic or when the
Gau{\ss} curvature of~$\fff$ is suitably non-degenerate, but the general
smooth case is still unsolved.

\subsection{Prescribed second fundamental form}\label{ssec: cartan}
In a different direction, in 1943 \'Elie Cartan studied the problem
of realizing a given second fundamental form~\cite{Ca:II}.  
In other words, he
studied the equation~$\sff_\xb = \sff$, where~$\sff$ is a given
quadratic form.  He showed that, when~$\sff$ is real-analytic 
and non-degenerate, the equation~$\sff_\xb=\sff$ is always locally
solvable.  Little seems to be known about this problem in the smooth
category or in the global setting.  Possibly this is because, 
as Cartan showed, this problem is never elliptic and, in fact, has rather
complicated characteristics.

\subsection{Bonnet surfaces}\label{ssec: bonnet surfaces}
One can imagine specifying other aspects of the data contained
in~$(\fff,\sff)$.  For example, if~$\kappa_1$ and~$\kappa_2$
are the eigenvalues of~$\sff$ with respect to~$\fff$, one can
imagine trying to find an~$\xb$ that realizes
a given~$(\fff,\kappa_1,\kappa_2)$.

This problem was first studied by Bonnet and then several other authors.
Of course, since the Gau{\ss} equation asserts that~$K=\kappa_1\kappa_2$
where~$K$ is the Gau{\ss} curvature of~$\fff$, this is an obvious necessary
condition for realizability, so suppose that this holds.  It turns out
that, even with this condition, the generic data~$(\fff,\kappa_1,\kappa_2)$
cannot be realized by an immersion.  This should be expected, since,
even with the Gau{\ss} equation restriction, the given data depends
essentially on four arbitrary functions of two variables, so some sort
of compatibility condition is necessary.

The most thorough local analysis was done by Cartan~\cite{Ca:Bonnet}, 
who showed that, for the generic data~$(\fff,\kappa_1,\kappa_2)$ 
(satisfying the Gau{\ss} equation) that does admit a normally oriented 
realization~$(\xb,\nb)$, such a realization is unique up to isometry.  
This uniqueness fails for three special classes of data:

First, there exists a special class of data~$(\fff,\kappa_1,\kappa_2)$, 
depending on four arbitrary functions of one variable, for which 
there exist exactly two normally oriented realizations~$(\xb_\pm,\nb_\pm)$
that are not Euclidean congruent.  In the recent literature,
these are called \emph{Bonnet pairs}~\cite{KPP}.  

For the second and third classes of data~$(\fff,\kappa_1,\kappa_2)$
to be described below, there exists a $1$-parameter family 
of normally oriented realizations~$(\xb_\theta,\nb_\theta)$.  

The second class consists of the data that are realizable by
surfaces of constant mean curvature.  In this case, the $1$-parameter
family containing a given immersion of constant mean curvature is
just the classical circle of associated surfaces (most well-known 
in the case of mean curvature zero, i.e., the minimal surfaces).
This class of data depends locally on two arbitrary functions of
one variable.

The third class consists of the data realizable by a $6$-parameter 
family of surfaces now known as the \emph{Bonnet surfaces}.  These
are not as easy to describe geometrically, so the reader is referred
to sources in the bibliography for further information, particularly
Cartan's article~\cite{Ca:Bonnet}, Chern's article~\cite{Ch:Bonnetsurf},
and the more recent article~\cite{BoEi}, where the relationship between
these surfaces and Painlev\'e equations is explored.

\subsection{Cartan's case studies}\label{ssec: Cartans case studies}
In fact, there are a large number of possible problems one could 
study about the existence and uniqueness of normally oriented 
realizations of partial data drawn from the first and second fundamental
forms.  In his famous 1945 memoir%
\footnote{In spite of the publication date, the reader might want 
to note that this memoir is actually based on the lecture notes 
of a course that Cartan gave in 1936-7 at the Facult\'e des
Sciences de l'Universit\'e de Paris. These lectures were attended
by S.-S. Chern, who was, at that time, a postdoctoral student.} 
\emph{Les syst\`emes diff\'erentiels ext\'erieurs et leurs 
applications g\'eometriques}, Cartan considered
a number of these problems as illustrations of his methods.  
In particular, Chapitre~VII is devoted to such problems
and is still one of the best sources for information about them.

\subsection{Prescribed shape operator}\label{ssec: presecribed shape op}
On particularly natural object one can construct from the 
data~$(\fff,\sff)$ is the \emph{Weingarten shape operator}. 
This is the linear mapping~$S:TD\to TD$ defined by
the relation
\begin{equation}
\sff(v,w) = \fff(v,S w) = \fff(S v, w).
\end{equation}
for any pair of tangent vectors~$v,w\in T_pD$.  Since~$S$ is
$\fff$-self-adjoint, it has real eigenvalues (which are, of
course, the principal curvatures of any normally oriented realization)
and is (pointwise) diagonalizable.  There are no other pointwise
conditions on~$S$.

Conversely, given an endomorphism~$S:TD\to TD$ of the tangent
bundle that is pointwise diagonalizable, one can consider the problem of 
finding a normally oriented immersion~$(\xb,\nb)$ whose shape operator 
is~$S$.  Since the choice of~$S$ is tantamount to choosing a section
of a bundle of rank~$4$ over~$D$, namely~$\End(TD)$, a shape operator
essentially depends on four arbitrary functions of two variables.
Thus, one does not expect to be able to realize every possible~$S$
as a shape operator.

For example, if~$S$ has equal eigenvalues at every point, so 
that~$S = \kappa\,\id_{TD}$ for some function~$\kappa$ on~$D$,
then~$S$ cannot be realized unless~$\kappa$ is a constant, since
the only totally umbilic surfaces in~$\E3$ are planes and spheres.
On the other hand, if~$\kappa$ is constant, then~$S$ is 
realized by a normally oriented immersion of~$D$ into a plane or
sphere of appropriate radius.  Thus, this case is trivial.

\subsubsection{Umbilics and rectangularity}\label{sssec: umbilic and rect}
It is natural to define the points of~$D$ at which~$S$ has two
equal eigenvalues to be the \emph{umbilic points} of~$S$.  The
presence of these points complicates the discussion, so, for
simplicity, I will assume that there are no $S$-umbilic points.
In this case, there will be two functions~$A>B$ on~$D$ so that
$A(p)$ and~$B(p)$ are eigenfunctions of~$S_p:T_pD\to T_pD$ and,
since~$D$ is simply connected, there will exist two vector fields
$\ab$ and~$\bb$ on~$D$ with dual $1$-forms~$\alpha$ and~$\beta$
so that
\begin{equation}\label{eq: S diagonalized}
S = A\,\,\ab\ot \alpha + B\,\,\bb\ot\beta\,. 
\end{equation} 
(Of course~$\ab$ and~$\bb$ are not unique.)

Locally, there exist coordinates~$(x,y)$ on~$D$
in which~$S$ has the more specific form
\begin{equation}\label{eq: S principal coords}
S = A(x,y)\,\,\frac{\partial\hfil}{\partial x}\otimes \d x
   +B(x,y)\,\,\frac{\partial\hfil}{\partial y}\otimes \d y\,.
\end{equation}
These two coordinates, the so-called \emph{$S$-principal coordinates}, 
are each unique up to reparametrization.  If $S$-principal
coordinates~$(x,y)$ can be chosen globally on~$D$ in such a
way that~$(x,y):D\to\R2$ embeds~$D$ as a coordinate rectangle in
the~$xy$-plane, then the pair~$(D,S)$ will be said to be \emph{rectangular}.

Since the study conducted in this article will be almost entirely
a local one, it does no harm to restrict to the umbilic-free, 
rectangular case, so this will often be assumed unless it is
specifically stated otherwise.

\begin{remark}[Computability 1]\label{rem: computability 1}
Given an endomorphism~$S:TD\to TD$, its eigenvalues and eigendirections
can be computed algebraically, so that, when~$S$ has distinct eigenvalues,
the form~\eqref{eq: S diagonalized} can be computed effectively.
However, finding principal coordinates~$(x,y)$ explicitly when one is
given an operator~$S$ in the form~\eqref{eq: S diagonalized}
requires one to solve two coupled, nonlinear ordinary differential equations,
something that \emph{cannot} be done effectively unless~$S$ has special
properties.

However, as will be seen, the computations that need to be done can
be done without the use of principal coordinates; they are merely
a convenient expository device.  The replacement, as will be seen,
is to use the eigenform decomposition of the $1$-forms that~$S$ induces:
Any $1$-form~$\phi$ on~$D$ can be written uniquely in the form~$\phi=\phi'
+\phi''$ where $\phi'$ is a multiple of~$\alpha$ and~$\phi''$ is a
multiple of~$\beta$.  Correspondingly, there is a decomposition of
the exterior derivative on functions:  $\d f = \d' f + \d'' f$ where 
$\d' f = (\d f)'$ and~$\d'' f = (\d f)''$.  Of course, these two 
operators can be computed algebraically from~$S$, without recourse
to differential equations.
\end{remark}

\subsubsection{Cartan's non-uniqueness analysis}
\label{sssec: Cartan non-unique}
In Probl\`eme~IX of Chapitre~VII of~\cite{Ca:EDS}, Cartan considers 
the generality of pairs of immersions~$\xb,\yb:D\to\E3$ that are noncongruent
but induce the same shape operator.  His analysis will be only summarized
here.  He shows that, modulo reparametrization, these pairs 
depend on six arbitrary functions of one variable.  From this, 
he concludes that the `generic' normally oriented immersion~$(\xb,\nb)$ 
is uniquely characterized up to Euclidean congruence by its shape operator.

Recently, Ferapontov~\cite{Fe:soliton1} has studied this non-uniqueness
problem from the point of view of integrable systems and has shown
that this problem (with some extra genericity
hypothesis, to be described more fully beow 
in Remark~\ref{rem: non-Frobenius}) 
is susceptible to being formulated 
as a Lax pair with a spectral parameter.

Cartan also shows that, if~$\xb:D\to\E3$ is an immersion that is free
of umbilics and has the property that one of its families of
principal curves is planar (as is the case, for example, 
for surfaces of revolution and, more generally, 
for the so-called molding surfaces), then the space of 
immersions~$\yb:D\to\E3$ that induce the same shape operator as~$\xb$
depends on one arbitrary function of one variable.  

The condition of having one of the families of principal curves be
planar is equivalent, in the local coordinate 
form~\eqref{eq: S principal coords} of the shape operator~$S$, 
to having either~$A_y=0$ or~$B_x=0$, i.e., one of the principal curvatures
should be constant along the orthogonal family of principal curves.

Cartan does not mention the 1933 
work of Finikoff and Gambier~\cite{Fi,FiGa}, 
and perhaps he was unaware of it.  Their work makes the same
observations about the shape operators of surfaces with one family
of principal curves being geodesics and they provide examples
of shape operators that can be realized in a $3$-parameter
family of distinct ways.  (They believed that they had
a classification of such, but, as Ferapontov points out 
in~\cite{Fe:classify3}, they missed an entire family,  
the one designated as Type~I in this article.)

Other recent results on uniqueness and non-uniqueness
for the prescribed shape operator problem can be found in~\cite{SVVW}
and~\cite{SVWW}.  The authors particularly study the case of 
surfaces of revolution and give examples that exhibit 
the non-uniqueness that shows up in Cartan's analysis.

Ferapontov's article~\cite{Fe:survey} is a valuable source of
information about the history of this problem, 
so the reader is referred there for more details.  

\subsubsection{New results}\label{sssec: new results}
In this article, the non-uniqueness problem will be examined
in detail and some new results will be proved about the
explicit computablity of shape operators with the maximum
degree of flexibility in their realizations.  The reader is 
reminded that the data~$(D,S)$ is assumed to be umbilic-free, 
rectangular, and smooth.

First, there is the observation 
(see Proposition~\ref{prop: affine structure}) 
that the space of congruence classes of
normally oriented immersions~$(\xb,\nb):D\to\E3$ that realize~$S$
has a natural affine structure in the sense that, if~$(\xb_0,\nb_0)$
and~$(\xb_1,\nb_1)$ both realize~$S$, then there is a naturally
constructed family~$(\xb_t,\nb_t)$ for~$0\le t\le 1$ of normally
oriented immersions defined up to Euclidean congruence that interpolates
between the two given immersions, and every immersion in this family
realizes~$S$ as its shape operator.  Moreover, if the two given
immersions are not Euclidean congruent, then any two distinct members 
of the family~$(\xb_t,\nb_t)$ are mutually incongruent.  This
affine structure was implicit already in the works of Finikoff and Gambier.

Second, in the case where~$A_y$ and~$B_x$ are nonvanishing on~$D$
(which is a generic condition), it turns out 
(see Theorem~\ref{thm: at most 3}) that the space
of Euclidean congruence classes of normally oriented 
immersions~$(\xb,\nb):D\to\E3$ realizing~$S$ can be naturally, 
affinely embedded as a convex set~$\cX(S)$ in an affine space of 
dimension $3$.  Moreover, this convex set~$\cX(S)$ will have an interior
if and only if~$A$ and~$B$ satisfy a system~$\cE(A,B)=0$ of four highly 
nonlinear partial differential equations, two of order three and and two
of order four.  

It is not clear \emph{a priori} that the overdetermined system~$\cE(A,B)=0$ 
has any solutions for which~$A_y$ and~$B_x$ are nonvanishing. 
Moreover, it is not difficult to show that this system is not involutive 
in Cartan's sense, so further analysis is needed to understand 
the local and global solutions.

Again, this system was implicit in the work of Finikoff and Gambier,
who first derived this upper bound on the dimension of the space of
realizations of a given shape operator.  However, it appears that their
analysis of it was flawed, as they missed a family of solutions.

If one makes the additional
assumption that~$A$ and~$B$ themselves are nonvanishing,%
\footnote{Of course, assuming that~$A_y$ and~$B_x$ are nonvanishing already
implies that the locus in~$D$ where either~$A$ or~$B$ vanishes
has no interior, so this is not a drastic assumption.}
then it is possible to reformulate the system~$\cE(A,B)=0$ 
in more geometric terms, so that its analysis becomes greatly simplified.
The key, already noticed by Finikoff and Gambier, is to deal with the
reciprocals~$U=1/A$ and~$V=1/B$ and then to define the 
coframing~$\theta=(\theta_1,\theta_2)$, where
\begin{equation}\label{eq: thetas in UV intro}
\theta_1 = \frac{\d' V}{V-U}\,,\qquad\qquad\qquad\qquad
\theta_2 = \frac{\d''U}{U-V}\,.
\end{equation} 
The system~$\cE(A,B)=0$ turns out to be equivalent to a lower
order (and much simpler) overdetermined system~$\cE'(\theta)=0$ 
for the coframing~$\theta$.  One then finds 
(see Theorem~\ref{thm: two types})
that the system~$\cE'(\theta)=0$ is satisfied
if and only if~$\theta$ satisfies 
one of three possible determined, involutive systems.  Once
the solutions to~$\cE'(\theta)=0$ are described, 
the equations~\eqref{eq: thetas in UV intro} can be regarded
as a first order, linear hyperbolic determined system 
for two functions~$U$ and~$V$ and standard techniques can then
be applied for its solution.

Using this simplification, one sees that
(see Theorem~\ref{thm: two types})
that the system~$\cE(A,B)=0$ is satisfied for~$AB\not=0$ 
if and only if~$A$ and~$B$ satisfy one of three possible 
determined, involutive systems.  

The first two of these systems are exchanged by the operation 
of exchanging~$(x,y)$ and~$(A,B)$, so they can be regarded 
as essentially equivalent.  Each of these systems consists 
of a second order equation and a third order equation.  Operators~$S$
that satisfy either of these systems will be referred to as being
of Type~I.  This is the type that was missed by Finikoff and Gambier
and first discovered by Ferapontov~\cite{Fe:classify3}.

The third system is invariant under the exchange of~$(x,y)$ and~$(A,B)$
and consists of a pair of second order equations that forms a hyperbolic
system for~$A$ and~$B$ in principal coordinates. Operators~$S$
that satisfy this system will be referred to as being
of Type~II.  These are the shape operators that were first found by
Finikoff and Gambier.

Since each of the systems~\eqref{eq: first system} 
and~\eqref{eq: second system} is a determined, involutive system, 
local solvability is easy to demonstrate.   Thus, there are
many examples of shape operators of either type.  In fact,
several explicit examples are given in this article.

The remainder of the article is devoted to the analysis of the
geometric properties of these two types of shape operators, particularly
with an eye to the explicit computability of their realizations.
Some of these results will now be described.

In the case of shape operators of Type~I, 
it is shown (see~\S\ref{sssec: nat prin coords I}) 
that there are essentially canonical principal coordinates~$(x,y)$
on~$D$ in which the system that defines them
can be linearized and explicitly integrated by the method of Darboux 
(see~\S\ref{sssec: euler linear I}).  
Thus, these shape operators can be regarded as explicitly known.  
One finds that, modulo reparametrization, the shape operators
of this type depend on three arbitrary functions of one variable.
In fact, one can do much better than an integration by the method
of Darboux:  One can place the general shape operator solution in canonical
principal coordinate form using only quadrature and the solutions
of a single linear second order ODE.

Moreover, the structure equations show that if~$S$ is of Type~I 
and~$(\xb,\nb):D\to\E3$ is any realization, 
then the Gau{\ss} images of one family of principal curves
are arcs of a $1$-parameter family of spherical circles in~$S^2$
whose curve of centers lies on a geodesic. 
(see Proposition~\ref{prop: spherical circles}).  
This leads to an explicit integration (up to quadrature%
\footnote{
The term \emph{quadrature} in this article carries its classical meaning: 
Quadrature is the operation of finding a primitive for a given closed 
$1$-form, i.e., given a closed $1$-form~$\phi$ on a simply connected 
domain~$D$, quadrature constructs a function~$f$ on~$D$ 
so that~$\d f = \phi$.  Symbolically, this is written~$f=\int\phi$.  
The classical authors regarded quadrature as an elementary operation
and devoted much energy to finding ways to solve differential equations
that only involved algebraic operations and quadrature.}%
) of the differential equations that determine 
the realizations~$\xb$ of~$S$.  Again, this is much sharper than
merely being able to linearize the realization equations.

In fact, in the case that the spherical images of \emph{both} 
families of principal curves are arcs of circles, the computation 
of the (local) realizations of~$S$ is reduced to a sequence 
of algebraic operations and quadratures
in a manner analogous to the Weierstra{\ss} formula for minimal surfaces
(\S\ref{sssec: weierstrass I}).

In the case of shape operators of Type~II, 
it is shown (see~\S\ref{sssec: nat prin coords II}) that
there are essentially canonical principal coordinates~$(x,y)$
on~$D$ in which the system that defines them can
be linearized.  In this case, the resulting linear system is not integrable
by the method of Darboux, though a representation due to Poisson can
be invoked to express the general solution, which depends on
two arbitrary functions of one variable.

What is particularly
remarkable is that the structure equations show that if~$S$ is of Type~II 
and~$(\xb,\nb):D\to\E3$ is any realization, then the Gau{\ss} image of 
the net of principal curves is a net of confocal spherical
ellipses in~$S^2$ (see Proposition~\ref{prop: Gauss prin curves II}).  
(This was known already to Finikoff and Gambier.)
Again, this allows one to reduce
the explicit integration of the differential equations that 
determine the (local) realizations~$\xb$ to a sequence 
of algebraic operations and quadratures analogous 
to the Weierstra{\ss} formula 
for minimal surfaces~(\S\ref{sssec: weierstrass II}).

Various examples of each Type are introduced and studied.  
Here are some highlights:
\begin{itemize}

\item An example~(see Example~\ref{ex: index zero type I}) 
is given of a surface of Type~I that has an isolated umbilic of index~0, 
thus showing that shape operator flexibility does not control 
the index of isolated umbilics.

\item The shape operators of either type that can
be realized by minimal immersions are determined 
and the corresponding surfaces are described 
(see Example~\ref{ex: minimal surface I}
and Example~\ref{ex: minimal surface II}).

\item The quadric surfaces with distinct principal axes 
belong to Type~II (see Example~\ref{ex: quadrics II}).

\item Complete examples of surfaces of either Type are 
given and compact convex examples are given%
\footnote{
This does \emph{not} mean that the surfaces are globally flexible
keeping the shape operator fixed, but only that each point in
the surface has a neighborhood that is flexible keeping the shape
operator fixed.  This may seem paradoxical, but it well illustrates
the different natures of the local and global problems.} 
of surfaces of Type~II (see Example~\ref{ex: linear type I} 
and Example~\ref{ex: polynomial II}). 

\end{itemize}

\subsection{Acknowledgements}\label{ssec: acknowledgements}
I want to thank Udo Simon for discussing his work and the
work of Martin Wiehe on the problem of prescribed shape operators.
It was those discussions that inspired the present article.

I must also thank Eugene Ferapontov, who read an earlier version
of this article and pointed out that I had reproduced many of the
results of Finikoff, Gambier and himself in my analysis.  I am
very grateful to him for supplying me with references and giving
me the chance to write a revised article in which proper historical 
credit is given.

I would also like to thank Editorial 
Board of \emph{Results in Mathematics} for the opportunity 
to contribute to a volume honoring Shiing-Shen Chern,
whose beautiful works on classical surface theory 
(and every branch of modern differential geometry as well) 
have inspired me throughout my career as a geometer.
I offer this article, whose topic and outlook are inspired by
Professor Chern's wonderful article~\cite{Ch:Bonnetsurf}, 
as a small token of my gratitude for his profound effect 
on my mathematical life.

\section[Differential Analysis]{The Differential Analysis}
\label{sec: diff anal}

The computations below will proceed by the method of the moving
frame, so the basic notation will be introduced here, along with
a few useful facts.

\subsection{The structure equations}\label{ssec: structure equations}
Let~$D$ be a simply connected surface, let~$\xb:D\to\E3$ 
be an immersion, and let~$\nb:D\to S^2$ be a choice of unit
normal, i.e., $\nb\cdot\nb=1$ and~$\nb\cdot\d\xb = 0$.  The
first and second fundamental forms are defined as before by
\begin{equation}
\fff = \d\xb\cdot\d\xb\,,\qquad\qquad\qquad
\sff = -\d\nb\cdot\d\xb\,.
\end{equation}
The immersion will be assumed to be free of umbilics, i.e., that~$\sff$
has two distinct eigenvalues with respect to~$\fff$ at every point. These
eigenvalues (the principal curvatures) will be denoted~$A$ and~$B$ and
it will be supposed that~$A>B$ throughout~$D$.  

There exist~$1$-forms~$\omega_1$ and~$\omega_2$ on~$D$ 
(unique up to a sign) so that~$\fff$ and~$\sff$ are diagonalized as
\begin{equation}
\fff = {\omega_1}^2+{\omega_2}^2\,,\qquad\qquad\qquad
\sff = A\,\,{\omega_1}^2+B\,\,{\omega_2}^2\,.
\end{equation}
The functions~$A$ and~$B$ are the principal curvatures.  The
corresponding shape operator~$S$ is given by the formula
\begin{equation}
S = A\,\ub_1\ot\omega_1 + B\,\ub_2\ot\omega_2
\end{equation}
where~$\ub_1$ and~$\ub_2$ are the vector fields on~$D$ dual
to the coframe field~$(\omega_1,\omega_2)$. (Note that
the sign ambiguity in the choice of~$\omega_1$ and~$\omega_2$
does not affect~$S$.)

Once the forms~$\omega_1$ and~$\omega_2$ are chosen, 
there will exist unique smooth mappings~$\eb_1,\eb_2:D\to S^2$
so that
\begin{equation}\label{eq: first str eq 1}
\d\xb = \eb_1\,\omega_1 + \eb_2\,\omega_2\,.
\end{equation}
The integral curves of the equation~$\omega_2=0$ map to tangents
to~$\eb_1$ and are called the first family of principal
curves while the integral curves of the equation~$\omega_1=0$ map 
to tangents to~$\eb_2$ and are called the second family of principal
curves.  The pair of foliations of~$D$ by the principal curves is
called the \emph{net} of principal curves induced by~$\xb$.

Setting~$\eb_3 = \nb$, the frame field~$(\eb_1,\eb_2,\eb_3)$ is
orthonormal, so the $1$-forms~$\omega_{ij} = \eb_i\cdot \d\eb_j$
satisfy~$\omega_{ij} = -\omega_{ji}$
and the equations
\begin{equation}\label{eq: first str eq 2}
\d\eb_i = \eb_1\,\omega_{1i} + \eb_2\,\omega_{2i} + \eb_3\,\omega_{3i}\,,
\qquad\qquad i = 1,2,3.
\end{equation}
They also satisfy
\begin{equation}\label{eq: sff principal components}
\omega_{31} = A\,\omega_1\,,\qquad\qquad\qquad
\omega_{32} = B\,\omega_2\,,
\end{equation}
and there is a relation of the form
\begin{equation}\label{eq: w12 formula}
\omega_{12} = u\,\omega_1 + v\,\omega_2
\end{equation}
for some functions~$u$ and~$v$ on~$D$.
In fact, at any point~$p\in D$, $u(p)$ is the geodesic
curvature of the first principal curve passing through~$p$
while~$v(p)$ is the geodesic curvature of the second principal
curve passing through~$p$.

These forms satisfy the structure equations
\begin{equation}\label{eq: second str eqs}
\begin{split}
\d\omega_1 & = {} -  \omega_{12} \w \omega_2\,,\\
\d\omega_2 & = \phantom{{} - {}} \omega_{12} \w \omega_1\,,\\[3pt]
\d\omega_{31} & = {} -  \omega_{12} \w \omega_{32}\,,\\
\d\omega_{32} & = \phantom{{} - {}}\omega_{12} \w \omega_{31}\,,\\[3pt]
\d\omega_{12} & = \phantom{{} - {}}\omega_{31} \w \omega_{32}\,.
\end{split}
\end{equation}

Conversely, the essential content of Bonnet's theorem is that, if~$\omega_1$,
$\omega_2$, $\omega_{31}$, $\omega_{32}$, and~$\omega_{12}$ are $1$-forms
defined on a simply connected surface~$D$, satisfy 
the equations~\eqref{eq: second str eqs},
and the relation~$\omega_{31}\w\omega_1+\omega_{32}\w\omega_2=0$
(which is a consequence of~\eqref{eq: sff principal components}), 
then there exist mappings~$\xb: D\to\E3$ 
and~$(\eb_1,\eb_2,\eb_3):D\to\Or(3)$
so that the structure equations~\eqref{eq: first str eq 1}
and~\eqref{eq: first str eq 2} hold and that such mappings are
unique up to composition with an isometry of~$\E3$.  If, in addition,
$\omega_1\w\omega_2$ is nonvanishing on~$D$, then~$\xb$ is an immersion.

\subsection{A first look}\label{ssec: first look} 
Suppose given an open rectangular domain~$D$ in the $xy$-plane 
and a smooth candidate for a shape operator
\begin{equation}
S = A(x,y)\, \frac{\partial\hfil}{\partial x}\otimes \d x
   +B(x,y)\, \frac{\partial\hfil}{\partial y}\otimes \d y\,,
\end{equation}
where~$A>B$ throughout~$D$.  

If~$S$ is realized by a normally oriented immersion~$\xb:D\to\E{3}$, 
there will exist a principal orthonormal frame 
field~$\eb = (\eb_1,\eb_2,\eb_3):D\to\Or(3)$ 
so that the structure equations
\begin{subequations}
\begin{align}
\d\xb &= \eb_1\,\omega_1 + \eb_2\,\omega_2\,,
\label{eq: dx str eq}\\
\d\eb_i  &= \eb_j\,\omega_{ji}\,,
\label{eq: de str eq}
\end{align}
\end{subequations}
hold, where, for some functions~$a,b >0$ on~$D$, the structure
forms satisfy
\begin{equation}
\omega_1 = \frac{\d x}{\sqrt{a}}\,,\quad
\omega_2 = \frac{\d y}{\sqrt{b}}\,,\quad
\omega_{31} = \frac{A\,\d x}{\sqrt{a}}\,,\quad
\omega_{32} = \frac{B\,\d y}{\sqrt{b}}\,.
\end{equation}
(This way of parametrizing the possible structure forms, which
may seem odd at first glance, turns
out to lead to a \emph{linear} inhomogeneous system of equations 
for $a$ and~$b$.)

The structure equations~$\d\omega_1 = - \omega_{12}\w\omega_2$
and~$\d\omega_2 = \omega_{12}\w\omega_1$ imply that
\begin{equation}\label{eq: first conn formula}
\omega_{12} = -\frac{a_y\,\sqrt{b}}{\sqrt{a^3}}\,\d x
              +\frac{b_x\,\sqrt{a}}{\sqrt{b^3}}\,\d y\,.
\end{equation}
The structure equations $\d\omega_{31} = - \omega_{12}\w\omega_{32}$
and~$\d\omega_{32} = \omega_{12}\w\omega_{31}$ then imply two 
linear equations for~$a$ and~$b$:
\begin{equation}\label{eq: first ab eqs}
a_y = \frac{2\,A_y}{(A{-}B)}\,a\,,\qquad b_x = \frac{2\,B_x}{(B{-}A)}\,b\,,
\end{equation}
so that \eqref{eq: first conn formula} can be written in the form
\begin{equation}
\omega_{12} = -\frac{2 A_y\,\sqrt{b}}{(A-B)\,\sqrt{a}}\,\d x
              +\frac{2 B_x\,\sqrt{a}}{(B-A)\,\sqrt{b}}\,\d y\,.
\end{equation}
The final structure equation~$\d\omega_{12} = \omega_{31}\w\omega_{32}$
then yields a third (inhomogeneous) linear equation for~$a$ and~$b$:
\begin{equation}\label{eq: last ab eq}
\begin{split}
-2\,AB\,(A{-}B)^2
 &=\phantom{{} + } (A{-}B)B_x\,a_x 
      -2\,\left(B_xA_x{-}2\,{B_x}^{2}{+}(B{-}A)B_{xx}\right)\,a\\
&\phantom{= } \, + (B{-}A)A_y\,b_y
      \,-2\,\left(A_yB_y{-}2\,{A_y}^{2}{+}(A{-}B)A_{yy}\right)\,b\,.
\end{split}
\end{equation}

For general functions~$A$ and~$B$, the equations \eqref{eq: first ab eqs} 
and \eqref{eq: last ab eq} define a system of three equations
for the two unknown functions~$a$ and~$b$ that is incompatible, 
i.e., there will be no solutions.  

Still, because
these equations are linear, inhomogeneous equations, it follows
that if~$(a_0,b_0)$ and~$(a_1, b_1)$ are positive solutions
to \eqref{eq: first ab eqs} and \eqref{eq: last ab eq}, then
setting
\begin{equation}
(a_t,b_t) = \left(\,(1-t)a_0+ta_1,\,(1-t)b_0+tb_1\,\right)
\end{equation}
for~$0\le t\le 1$ defines a `segment' of positive solutions.
Since~$D$ is simply connected, the following
result is a direct consequence of Bonnet's theorem.

\begin{proposition}\label{prop: affine structure}
The space of Euclidean congruence classes of normally oriented 
immersions~$\xb:D\to\E{3}$ that realize the shape operator~$S$ 
is a convex set in the affine space consisting of the solutions 
of the inhomogeneous linear system defined by~\eqref{eq: first ab eqs} 
and~\eqref{eq: last ab eq}. \qed
\end{proposition}

There remains the question of determining conditions on~$A$ and~$B$
that will determine whether or not there exist any solutions
to the system \eqref{eq: first ab eqs} and \eqref{eq: last ab eq}. 
A full analysis of their compatibility has to be broken into a number
of cases.  As will be seen, determining necessary and sufficient 
conditions for the existence of a single solution is likely to be
rather complicated.

\subsection{Two elementary cases}\label{ssec: elementary cases}
In the first place, if $A_y=B_x \equiv 0$ then \eqref{eq: last ab eq} is
incompatible unless one of $A$ or~$B$ also vanishes.  By symmetry,
one can assume that~$A\equiv 0$.   Then \eqref{eq: last ab eq} is
an identity and the relations~\eqref{eq: first ab eqs} are equivalent
to the conditions
\begin{equation}
a = u(x)\qquad b = v(y)
\end{equation}
for some positive functions~$u$ and~$v$ of a single variable.  
The corresponding surfaces are generalized cylinders and no further 
discussion is required.

In the second place, one could have exactly one of~$A_y\equiv0$ 
or~$B_x\equiv0$.  Again, by symmetry it suffices to treat the
case~$A_y\equiv0$ under the assumption that~$B_x\not=0$.  Note that,
geometrically, this condition corresponds to the case where the
(principal) $\eb_1$-curves are congruent planar geodesics.  The
discussion to be given below is essentially due to Cartan, who used
somewhat different terminology and notation.

In this case, one must have~$a = u(x)$ for some positive function~$u$ 
of one variable, just as before, but now \eqref{eq: last ab eq} is not
an identity; instead, it becomes the inhomogeneous linear equation
\begin{equation}\label{eq: last ab eq A_y=0}
0 = u'(x)
+2\,\frac{\left(B_xA_x{-}2\,{B_x}^{2}{+}B_{xx}(B{-}A)\right)}
              {(B{-}A)B_x}\,u(x)
-2\,\frac{AB\,(A{-}B)^2}{(B{-}A)B_x}\,,
\end{equation}
which, for notational simplicity, can be written in the form
\begin{equation}\label{eq: last ab eq A_y=0 simplified}
0 = u'(x) + f(x,y)\,u(x) + g(x,y).
\end{equation}
Any solution~$u$ to~\eqref{eq: last ab eq A_y=0 simplified}
must also satisfy its derivative with respect to~$y$:
\begin{equation}\label{eq: last ab eq A_y=0 simplified dy}
0 = f_y(x,y)\,u(x)+g_y(x,y).
\end{equation}

If~$f_y(x,y)\equiv0$ but~$g_y(x,y)\not=0$, then there is
no solution to~\eqref{eq: last ab eq A_y=0 simplified dy} 
and hence no realization of~$S$ as a shape operator.

If~$f_y(x,y)$ is non-zero, there can be no more than one solution
to \eqref{eq: last ab eq A_y=0 simplified dy} and this may or
may not be positive and may or may not satisfy 
\eqref{eq: last ab eq A_y=0 simplified}.  If there is no solution,
then~$S$ cannot be realized as a shape operator.  If there is
a solution and it is either not positive or does not 
satisfy~\eqref{eq: last ab eq A_y=0 simplified}, 
then there is no positive solution to~\eqref{eq: last ab eq A_y=0 simplified}
and hence no realization of~$S$.  

On the other hand, if \eqref{eq: last ab eq A_y=0 simplified dy} is
an identity, then~$f$ and $g$ depend only on~$x$, so that
\eqref{eq: last ab eq A_y=0 simplified} is an ordinary differential
equation for~$u$ that has a $1$-parameter family of solutions.%
\footnote{
This happens, for example, in the case of surfaces of revolution,
where both~$A$ and~$B$ are functions of~$x$ alone.  As Cartan
points out, it also happens for the more general case of \emph{molding
surfaces}.}   In particular, there must be positive solutions,
at least in open $x$-intervals.

In any case, if~\eqref{eq: last ab eq A_y=0} has a positive solution~$u$,  
the remaining equation to be satisfied is the homogeneous linear equation 
for~$b$
\begin{equation}
b_x = \frac{2\,B_x}{(B-A)}\,b\,,
\end{equation}
whose general positive solution is of the form~$b(x,y) = v(y)\,\bar b(x,y)$ 
where~$\bar b>0$ is any particular solution and $v$ is an arbitrary positive 
function of one variable.

Thus, the positive solutions of~\eqref{eq: first ab eqs} 
and~\eqref{eq: last ab eq} (when they exist)
are seen to depend essentially on one arbitrary function of one
variable (plus possibly one constant).  It is in this sense that
Cartan means his statement that the realizations of~$S$ in such cases
depend on one arbitrary function of one variable.

For more information about these cases, along with an interesting
discussion of examples, the reader can consult~\cite{SVVW} and~\cite{SVWW}.

\subsection{The nondegenerate case}\label{ssec: nondegenerate case}
Now consider the general case, where~$A_y$ and~$B_x$ are nonvanishing.
In the notation of Remark~\ref{rem: computability 1}, 
this is equivalent to the assumption that~$\d''A$ and~$\d'B$ are
nonvanishing (and hence is checkable without having to find
principal coordinates beforehand).  Geometrically, this is equivalent
to the condition that the principal curves of any realization~$\xb:D\to\E3$
should have nonvanishing geodesic curvature.

The relation \eqref{eq: last ab eq} can now be expressed in the form
\begin{equation}\label{eq: last ab eq in p}
\begin{split}
a_x& = \left(\frac{1}{B_x}\right)\,p
      -2\left({\frac {B_{xx}}{B_x}}+{\frac{2\,B_x{-}A_x}{(A{-}B)}}\right)\,a
      +{\frac {AB(B{-}A)}{B_x}}\,,\\
b_y& = \left(\frac{1}{A_y}\right)\,p
      -2\left ({\frac{A_{yy}}{A_y}}+{\frac{2\,A_y{-}B_y}{(B{-}A)}}\right)\,b
      +{\frac {AB(A{-}B)}{A_y}}\,,
\end{split}
\end{equation}
for some unknown function~$p$.  

Differentiating the first equation of~\eqref{eq: first ab eqs} with
respect to~$x$ and the first equation of~\eqref{eq: last ab eq in p}
with respect to~$y$ and comparing the two expressions for~$a_{xy}$ leads
to an equation of the form
\begin{equation}\label{eq: p_y formula}
p_y = \frac{E_1\,p + E_2\,a + E_3}{A_y(A-B)^2}
\end{equation}
where~$E_1$,~$E_2$, and~$E_3$ are certain polynomials in the 
terms
$$
A,A_x,A_y,B,B_x,B_y,B_{xx},B_{xy},B_{xxy}
$$ 
whose exact form will not be needed in the discussion below.  Similarly,
differentiating the second equation of~\eqref{eq: first ab eqs} with
respect to~$y$ and the second equation of~\eqref{eq: last ab eq in p}
with respect to~$x$ and comparing the two expressions for~$b_{xy}$ leads
to an equation of the form
\begin{equation}\label{eq: p_x formula}
p_x = \frac{E_4\,p + E_5\,b + E_6}{B_x(A-B)^2}
\end{equation}
where~$E_4$,~$E_5$, and~$E_6$ are polynomials in the 
terms
$$
B,B_x,B_y,A,A_x,A_y,A_{yy},A_{xy},A_{xyy}\,.
$$

Considering the combined equations~\eqref{eq: first ab eqs}, 
\eqref{eq: last ab eq in p}, \eqref{eq: p_y formula}, and  
\eqref{eq: p_x formula} as a total system of equations for the
three unknowns~$a$, $b$, and $p$, one sees that there is at
most a three parameter family of solutions.  In fact, a solution,
if it exists, is uniquely determined by specifying the values 
of~$a$, $b$, and~$p$ at a single point~$(x_0,y_0)$ in~$D$. 
This, combined with Bonnet's theorem, yields the following 
fundamental result.

\begin{theorem}\label{thm: at most 3}
If~$S$ satisfies~$A_y\not=0$ and~$B_x\not=0$,
then the space of Euclidean congruence classes of normally oriented
immersions~$\xb:D\to\E3$ that realize~$S$ is a convex set 
in a vector space of dimension~$3$. \qed
\end{theorem}

For general~$A$ and~$B$, the combined system~\eqref{eq: first ab eqs}, 
\eqref{eq: last ab eq in p}, \eqref{eq: p_y formula}, and  
\eqref{eq: p_x formula} will not have any solutions at all.
This article is devoted to understanding the exceptional case 
in which this combined system is Frobenius, i.e., 
for which it possesses a $3$-parameter family of solutions.  

By the usual Frobenius criterion, this is the
case if and only if the two expressions for~$p_{xy}$ got
by differentiating \eqref{eq: p_y formula} with respect to~$y$ 
and by differentiating \eqref{eq: p_x formula} with respect to~$x$ 
agree (after, of course, taking into account the full system
of equations).  Carrying out this comparison~$(p_x)_y=(p_y)_x$ 
yields an equation of the form
\begin{equation}
\frac{E_7\,a + E_8\,b + E_9\,p + E_{10}}{{A_y}^2{B_x}^2(A-B)^3} = 0
\end{equation}
where $E_7$, $E_8$, $E_9$, and $E_{10}$ are polynomials in~$A$, $B$,
and certain of their derivatives up to and including order~$4$.
Thus, the Frobenius criterion is satisfied if and only if these
four polynomial expressions vanish identically:
\begin{equation}
E_7\equiv E_8\equiv E_9\equiv E_{10} \equiv 0.
\end{equation}
This is an overdetermined system~$\cE(A,B)=0$ for the functions~$A$ and~$B$.
In fact,~$E_9=0$ and~$E_{10}=0$ are third order equations for~$A$ and~$B$
while~$E_7$ and~$E_8$ are fourth order. These four
expressions are rather complicated.  For example, $E_7$ and $E_8$
have~$29$ terms apiece, $E_9$ has~$18$ terms, and~$E_{10}$ has~$42$ terms.

This system is  not involutive, so that determining 
its space of solutions requires further study.  
It is possible to study these equations directly, but it turns out
that the analysis is simplified and made more geometric by making
a change of variables and redoing the calculation up to this point.
That is the subject of the next subsection.

\subsection{A second look}\label{ssec: second look} 
Now suppose given a shape operator~$S$ of the 
form~\eqref{eq: S diagonalized} on a simply-connected surface~$D$
(that is not necessarily rectangular with respect to~$S$).

\subsubsection{Nondegeneracy}\label{sssec: nondegeneracy}
I suppose, as in the previous subsection, 
that~$A-B$ is nonvanishing on~$D$ and that~$\d''A$ 
and~$\d'B$ are also nonvanishing on~$D$.  An operator~$S:TD\to TD$
satisfying these conditions will be said to be \emph{nondegenerate}.

\subsubsection{Invertibility}\label{sssec: invertibility}
In addition to implying that~$\d''A$  and~$\d'B$ are everywhere
linearly independent on~$D$, these hypotheses imply 
that the zero locus of~$A$ (if non-empty) 
consists of smooth curves transverse to the second family 
of principal curves and that the zero locus of~$B$ (if non-empty)
consists of smooth curves transverse to the first family 
of principal curves.  In particular, 
the open set in~$D$ where~$AB\not=0$, i.e., 
where $S$ is invertible, is dense.

While the analysis below can be carried out in a neighborhood
of a curve in~$D$ along which~$A$ and/or $B$ vanishes, this considerably
complicates the discussion and does not seem to be worth the trouble.
Thus, for simplicity of exposition, I am going to further
impose the condition that~$A$ and~$B$ themselves be nonvanishing
on~$D$, i.e., that~$S$ be invertible on all of~$D$.

\subsubsection{A canonical coframing}\label{sssec: canonical coframing}
It will be useful to define two $1$-forms
\begin{equation}\label{eq: thetas defined}
\theta_1 =   \frac{A\,\d' B}{B(B-A)}\,,\qquad\qquad\qquad
\theta_2 =   \frac{B\,\d''A}{A(A-B)}\,.
\end{equation}
Note that these two $1$-forms constitute a coframing 
on~$D$ that depends only on~$S$ (and not on any choice
of principal coordinates).  In fact, $\theta_1$ and~$\theta_2$
depend on one derivative of~$S$.%
\footnote{More precisely, $\theta_1$ and~$\theta_2$ are algebraic
functions of the $1$-jet of~$S$.}

It also turns out to be more convenient to work 
with the reciprocals of~$A$ and~$B$ than with~$A$ and~$B$ themselves.
Thus, set
\begin{equation}\label{eq: UV defined}
U = \frac1A\,,\qquad\qquad V = \frac1B\,.
\end{equation}
In terms of~$U$ and~$V$, the $\theta_i$ have the expressions
\begin{equation}\label{eq: thetas in UV}
\theta_1 = \frac{\d' V}{V-U}\,,\qquad\qquad\qquad\qquad
\theta_2 = \frac{\d''U}{U-V}\,.
\end{equation}

Since these forms are linearly independent, 
there exist unique functions~$K_1$ and~$K_2$ on~$D$ 
so that
\begin{equation}\label{eq: K defined}
\d\theta_1 =   K_1\,\theta_1\w\theta_2\,,\qquad\qquad\qquad
\d\theta_2 =   K_2\,\theta_2\w\theta_1\,.
\end{equation}
The functions~$K_1$ and~$K_2$ depend on two derivatives of~$S$.

\subsubsection{Derivation of the total differential system}
\label{sssec: derivation}
Now, if~$S$ is to be induced by an immersion~$\xb:D\to\E{3}$, 
there will exist a principal orthonormal frame 
field~$\eb = (\eb_1,\eb_2,\eb_3):D\to\Or(3)$ so that the
structure equations~\eqref{eq: dx str eq} and~\eqref{eq: de str eq}
hold, where, for some positive functions~$a$~and~$b$ on~$D$, 
the structure forms satisfy%
\footnote{Obviously, these equations are got from the original
equations by `rescaling' the original~$a$ and~$b$.  This rescaling
greatly simplifies the calculations, a fact that
was only noticed in hindsight.}
\begin{equation}\label{eq: w1,w2,w31,w32 normalized}
\begin{split}
\omega_1 &= \frac{U\,\theta_1}{\sqrt{a}}\,,\qquad\qquad
\omega_{31} = \frac{\theta_1}{\sqrt{a}}\,,\\
\omega_2 &= \frac{V\,\theta_2}{\sqrt{b}}\,,\qquad\qquad
\omega_{32} = \frac{\theta_2}{\sqrt{b}}\,.
\end{split}
\end{equation}

Computing as in the previous sections, one finds that 
the structure equations for~$\d\omega_1$,~$\d\omega_2$,
$\d\omega_{31}$,~$\d\omega_{32}$, and~$\d\omega_{12}$ are
equivalent to the condition that 
\begin{equation}\label{eq: w12 normalized}
\omega_{12} = {} -\frac{\sqrt{b}}{\sqrt{a}}\,{\theta_1}
                 +\frac{\sqrt{a}}{\sqrt{b}}\,{\theta_2}
\end{equation}
where~$a$ and~$b$ satisfy the equations
\begin{equation}\label{eq: da db}
\begin{split}
\d a & =    \phantom{(2b+p+1)}\hbox to 0pt{\hss$(2a+p+1)$}\,\theta_1
         +  \phantom{(2a-p+1)}\hbox to 0pt{\hss$2a\,(1-K_1)$}\,\theta_2\,,\\
\d b & =    \phantom{(2b+p+1)}\hbox to 0pt{\hss$2b\,(1-K_2)$}\,\theta_1 
         +  \phantom{(2a-p+1)}\hbox to 0pt{\hss$(2b-p+1)$}\,\theta_2\,.\\
\end{split}
\end{equation}
for some function~$p$ on~$D$.  

Before these equations can be differentiated, it will be necessary
to introduce derivatives of the functions~$K_i$ in the form
\begin{equation}\label{eq: dK_i}
\d K_1 = K_{11}\,\theta_1+K_{12}\,\theta_2\,,\qquad\qquad
\d K_2 = K_{21}\,\theta_1+K_{22}\,\theta_2\,.
\end{equation}

Now, taking the exterior derivative of the equations~\eqref{eq: da db}
shows that the exterior derivative of~$p$ must be given by
\begin{equation}\label{eq: dp}
\begin{split}
\d p &= \bigl((2-K_2)(p-1)+2(K_1-K_2-K_1K_2+K_{22})\,b\bigr)\,\theta_1\\
     &\phantom{=} 
      + \bigl((2-K_1)(p+1)+2(K_1-K_2+K_1K_2-K_{11})\,a\bigr)\,\theta_2\,.
\end{split}
\end{equation}

\begin{remark}[An uncoupling]\label{rem: uncoupling}
The reader should note the very interesting fact that 
equations~\eqref{eq: da db} and~\eqref{eq: dp} do not involve~$U$
or~$V$ directly, but are expressed solely in terms of the 
coframing~$\theta$ and its derivatives.  
This is important for two reasons:  

First, the coframing~$\theta$ 
contains less information than the operator~$S$ and so has fewer
local invariants.  This is useful because the local compatibility
conditions for \eqref{eq: da db} and~\eqref{eq: dp} can now
be expressed in terms of the local invariants of the 
coframing~$\theta$, as will be seen below.

Second, because the system \eqref{eq: da db} and~\eqref{eq: dp}
is expressed in terms of the first and second derivatives of~$\theta$,
its compatibility conditions will be expressed in terms of third
derivatives of~$\theta$, which is much simpler than fourth derivatives
of~$S$.  Thus, it is a lower order problem in these terms.
\end{remark}

\subsubsection{The Frobenius condition}
\label{sssec: the Frobenius condition}
Introducing the functions~$K_{ijk}$ by the equations
\begin{equation}\label{eq: dK_ij}
\d K_{ij} = K_{ij1}\,\theta_1+K_{ij2}\,\theta_2\,,
\end{equation}
the exterior derivative of each side of~\eqref{eq: dp} 
can now be computed.  The result is that the following 
inhomogeneous linear relation among $a$, $b$, and~$p$ must hold
\begin{equation}\label{eq: abp linear relation}
\begin{split}
0&= 8 - 2\,K_1 - 2\,K_2 - 4\,K_1K_2 + 3\,K_{11}+3\,K_{22}\\
&\quad {}- 2\bigl((1{+}K_2)K_{11}-(1{-}K_1)K_{21}-K_{111}\bigr)\,a\\
&\quad {}+ 2\bigl((1{-}K_2)K_{12}-(1{+}K_1)K_{22}+K_{222}\bigr)\,b\\
&\quad {}+ 3(K_{11}-K_{22}-2\,K_1+2\,K_2)\,p\,.
\end{split}
\end{equation}

Thus, the necessary and sufficient condition on the operator~$S$ %
\footnote{more precisely, on the coframing~$\theta$} 
in order that the combined system given by~\eqref{eq: da db} 
and~\eqref{eq: dp} be Frobenius (and hence have a three parameter 
family of solutions) is that the following four equations hold:
\begin{equation}\label{eq: Frob conds}
\begin{split}
0 &= 8 - 2\,K_1 - 2\,K_2 - 4\,K_1K_2 + 3\,K_{11}+3\,K_{22}\,,\\
0 &= K_{11}-K_{22}-2\,K_1+2\,K_2\,,\\
0 &= (1{+}K_2)K_{11}-(1{-}K_1)K_{21}-K_{111}\,,\\
0 &= (1{-}K_2)K_{12}-(1{+}K_1)K_{22}+K_{222}\,.
\end{split}
\end{equation}
For the rest of this section, equations~\eqref{eq: Frob conds}
will be assumed.

\begin{remark}[Non-Frobenius cases]\label{rem: non-Frobenius}
Although this article will contain no further discussion of the 
non-Frobenius case, it might be helpful to the reader to have
some remarks about how the analysis can be continued to classify
the operators for which the space of solutions to \eqref{eq: da db} 
and~\eqref{eq: dp} has dimension less than~$3$.

In the first place, if the relation~\eqref{eq: abp linear relation}
is non-trivial, then the space of solutions~$(a,b,p)$ 
to~\eqref{eq: da db} and~\eqref{eq: dp} has dimension at most~$2$.

In the general case, where the coefficients of~$a$, $b$, and~$p$
in~\eqref{eq: abp linear relation} are not identically vanishing,
applying the exterior derivative to~\eqref{eq: abp linear relation}
and taking the coefficients of~$\theta_1$ and~$\theta_2$ on the
right hand side will yield two more relations of the form
\begin{equation}\label{eq: differentiated abp relations}
\begin{split}
0 &= A_1\,a + B_1\,b+ P_1\,p + Q_1\,, \\
0 &= A_2\,a + B_2\,b+ P_2\,p + Q_2\,.
\end{split}
\end{equation}
where the coefficients~$A_i$, $B_i$, $P_i$, and~$Q_i$ are polynomials
in~$K_1$, $K_2$ and their coframing derivatives up to order~$4$.  

The combined linear system~\eqref{eq: abp linear relation} 
and~\eqref{eq: differentiated abp relations} will generically have
a unique solution~$(a,b,p)$.  This solution may or may not have~$a$
and~$b$ positive on~$D$ and, even if they are positive, this
solution may not satisfy~\eqref{eq: da db} and~\eqref{eq: dp}.
If $a$ and~$b$ are positive and the system~$(a,b,p)$ does 
satisfy~\eqref{eq: da db} and~\eqref{eq: dp}, then, by Bonnet's
Theorem,~$S$ can be realized as a shape operator and in only one way.

However, it can happen that the combined 
equations~\eqref{eq: abp linear relation} 
and~\eqref{eq: differentiated abp relations} are a linearly
dependent system, admitting either a one-dimensional or two-dimensional
space of (algebraic) solutions.  

For example, it admits a two-dimensional space of solutions if and
only if the equations in~\eqref{eq: differentiated abp relations} are
multiples of the equation~\eqref{eq: abp linear relation}.  In this
case, the combined system~\eqref{eq: da db} and~\eqref{eq: dp} is
Frobenius when restricted to the relation~\eqref{eq: abp linear relation}.
As a result, the combined system~\eqref{eq: da db} and~\eqref{eq: dp}
has a two dimensional space of solutions.  In this case, any
point of~$D$ where there is a solution to~\eqref{eq: abp linear relation}
with~$a$ and~$b$ positive at the given point has a neighborhood
on which~$S$ can be realized as a shape operator and in a two-parameter
family of ways.  Note that the condition that the 
equations~\eqref{eq: differentiated abp relations} be
multiples of the equation~\eqref{eq: abp linear relation} is a system
of fourth order PDE for the coframing~$\theta$.  This system has now
been partially analyzed, but the description of its space of solutions
is complicated.  This will be the subject of a future article.

Finally, if the combined system~\eqref{eq: abp linear relation} 
and~\eqref{eq: differentiated abp relations} has rank one, then
a further differentiation will test whether the restriction 
of~\eqref{eq: da db} and~\eqref{eq: dp} to the solutions of these
relations is Frobenius.  If this is the case and there is a solution
with~$a$ and~$b$ positive, then~$S$ can be realized as a shape operator
in a one-parameter family of ways.  This condition is seen to be a
set of fifth order PDE for the coframing~$\theta$, but has not yet been
fully analyzed.  However, Cartan's analysis of the non-uniqueness
part of the problem mentioned in~\S\ref{sssec: Cartan non-unique}
can be applied in this case to show that the space of such operators 
must depend on four arbitrary functions of one variable. (The right
count is four, not six, because two arbitrary functions are lost
in the passage from~$S$ to the coframing~$\theta$.)  In fact,
Ferapontov~\cite{Fe:soliton1} has shown that this system can be
cast into the standard framework of an integrable system described
as a Lax pair with a parameter.  
\end{remark}

\subsubsection{Consequences}
\label{sssec: Frobenius consequences}
It is not at all clear how many nondegenerate invertible operators~$S$ 
there are that satisfy the four conditions~\eqref{eq: Frob conds}.  
The following analysis will show that these conditions imply one
of three possible determined systems.

The first two of the equations~\eqref{eq: Frob conds} 
do not involve the~$K_{ijk}$ and can be written in the form
\begin{equation}\label{eq: Frob conds 2}
K_{11} = {\ts\frac23}(K_1-1)(K_2+2)\,,\qquad\qquad
K_{22} = {\ts\frac23}(K_1+2)(K_2-1)\,.
\end{equation}
These two equations have an important consequence:  
The first can be
written in the form
\begin{equation}\label{eq: K_1 vanishes or not}
\d(K_1{-}1)\equiv{\ts\frac23}(K_2{+}2)(K_1{-}1)\,\theta_1\mod\theta_2\,.
\end{equation}
It thus follows, by the uniqueness of solutions of ordinary differential
equations, that~$K_1{-}1$ vanishes at a point of~$D$ if and only if it
vanishes along the entire first principal curve passing through that point.
Similarly, the second equation implies that~$K_2{-}1$ vanishes
at a point of~$D$ if and only if it vanishes along the entire second
principal curve passing through that point.

Using the equations~\eqref{eq: Frob conds 2} to eliminate~$K_{11}$
and~$K_{22}$, the formula~\eqref{eq: dp} for~$\d p$ can be simplified 
so that it does not involve any of the~$K_{ij}$. 
Using this new formula to recompute the identity~$\d(\d p) = 0$
then yields the relation
\begin{equation}\label{eq: abp linear relation 2}
0 = (K_1{-}1)(3K_{21}+2{K_2}^2+2K_2-4)\,a
   +(K_2{-}1)(3K_{12}+2{K_1}^2+2K_1-4)\,b\,,
\end{equation}
which must be an identity, i.e., the coefficients of~$a$ and~$b$
must vanish identically. Thus, equations~\eqref{eq: Frob conds} imply
the second order system formed by~\eqref{eq: Frob conds 2} and
\begin{equation}\label{eq: Frob conds 3}
\begin{split}
   (K_1-1)(3K_{21}+2{K_2}^2+2K_2-4) &= 0,\\
   (K_2-1)(3K_{12}+2{K_1}^2+2K_1-4) &= 0\,.
\end{split}
\end{equation}

The analysis must now be broken into a few separate cases.

\subsubsection{A genericity assumption}
\label{sssec: genericity assumption}
First, consider the open set~$D_0\subset D$ on which $(K_1-1)(K_2-1)$
is nonzero.  Then by~\eqref{eq: Frob conds 3}, the following
equations must hold on~$D_0$:
\begin{equation}\label{eq: Frob conds D_0}
\begin{split}
K_{12} = -{\ts\frac23}(K_1+2)(K_1-1)\,,\qquad\qquad
K_{21} = -{\ts\frac23}(K_2+2)(K_2-1)\,.
\end{split}
\end{equation}
Combined with~\eqref{eq: Frob conds 2}, this gives the formulae
\begin{equation}\label{eq: dK_i D_0}
\begin{split}
\d K_1 & = +{\ts\frac23}(K_1-1)(K_2+2)\,\theta_1 
           -{\ts\frac23}(K_1+2)(K_1-1)\,\theta_2\,,\\
\d K_2 & = -{\ts\frac23}(K_2+2)(K_2-1)\,\theta_1 
           +{\ts\frac23}(K_1+2)(K_2-1)\,\theta_2\,.\\
\end{split}
\end{equation}
Taking the exterior derivative of these equations yields the
relations
\begin{equation}\label{eq: D_0 final}
(K_1-1)(K_1+2)(K_2+2) = (K_2-1)(K_1+2)(K_2+2) = 0.
\end{equation}
Since~$(K_1-1)(K_2-1)$ is nonvanishing on~$D_0$, it
follows that~$(K_1+2)(K_2+2)$ vanishes identically.  

In fact, both~$(K_1+2)$ and~$(K_2+2)$ must vanish identically on~$D_0$.
For example, if~$(K_1+2)$ were nonzero at a point of~$D_0$, then~$(K_2+2)$
must vanish on a neighborhood of this point.  The second equation
of~\eqref{eq: dK_i D_0} then gives a contradiction since
the left hand side vanishes on this neighborhood but
the coefficient of~$\theta_2$ on the right hand side cannot be zero.
Similarly, $(K_2+2)$ must vanish at every point of~$D_0$.
Thus,~$K_1$ and~$K_2$ are each constant and equal to~$-2$ on~$D_0$.  
It then follows from the connectedness of~$D$ 
that either~$D_0$ is either empty or equal to all of~$D$.  

In summary, if~$D_0$ is nonempty, then~$D_0=D$, 
and~\eqref{eq: K defined} simplifies to
\begin{equation}\label{eq: second system}
\d\theta_1 =   -2\,\theta_1\w\theta_2\,,\qquad\qquad
\d\theta_2 =   -2\,\theta_2\w\theta_1\,.
\end{equation}
Moreover, the total differential system for~$a$, $b$, and~$p$ 
simplifies to the system
\begin{equation}\label{eq: da db dp on D_0}
\begin{split}
\d a & =    \phantom{{} - 4(2b-p+1)}\hbox to 0pt{\hss$(2a+p+1)$}\,\theta_1
         +  \phantom{4(2a+p+1)}\hbox to 0pt{\hss$6a$}\,\theta_2\,,\\
\d b & =    \phantom{{} - 4(2b-p+1)}\hbox to 0pt{\hss$6b$}\,\theta_1 
         +  \phantom{4(2a+p+1)}\hbox to 0pt{\hss$(2b-p+1)$}\,\theta_2\,,\\
\d p &=   {} - 4(2b-p+1)\,\theta_1 + 4(2a+p+1)\,\theta_2\,,
\end{split}
\end{equation}
which, in view of~\eqref{eq: second system}, is Frobenius.  

No further information can be gained by differentiating the
equations~\eqref{eq: second system} or~\eqref{eq: da db dp on D_0},
since these merely yield identities.  In \S\ref{ssec: Type II},
the systems of this type will be explicitly described.

\subsubsection{A special assumption}\label{sssec: special assumption}
On the other hand, suppose that~$D_0$ is empty, i.e., that 
$(K_1-1)(K_2-1)$ vanishes identically on~$D$.  Let~$R\subset D$
be any open principal rectangle in~$D$.   

Suppose that~$(K_1-1)$ 
does not vanish identically on~$R$.  Then, by the remark 
following~\eqref{eq: K_1 vanishes or not}, there
exits a principal curve in the first family on which~$(K_1-1)$ 
is nowhere vanishing.  Thus,~$(K_2-1)$ must vanish identically on 
this curve and so, again by the remark 
following~\eqref{eq: K_1 vanishes or not},~$(K_2-1)$
must vanish identically on~$R$.  (Since~$R$ is a principal
rectangle, every principal curve from the first family meets every
principal curve of the second family.)  

Similarly, if~$(K_2-1)$ does not vanish vanish identically on~$R$,
then~$(K_1-1)$ must vanish identically on~$R$.  

Since~$D$ is the union of its open principal rectangles,
it follows that~$D$ is the union of two open sets:
\begin{enumerate}
\item $D_1$,
on which~$(K_1-1)$ vanishes identically, and
\item $D_2$, on which~$(K_2-1)$
vanishes identically.%
\footnote{The reader should not jump to the conclusion that~$D_i$ is
equal to the locus where~$(K_i-1)$ vanishes.  It is only being asserted
that the \emph{interiors} of the two zero loci form a covering of~$D$.}
\end{enumerate}

Consider the open set~$D_1$.  Since~$(K_1-1)$ vanishes identically on~$D_1$,
it follows that~$K_{11}$ and~$K_{12}$ must also vanish identically on~$D_1$.
Thus, both equations in~\eqref{eq: Frob conds 3} are satisfied, 
so the $abp$ system is Frobenius.  

To summarize, on the open set~$D_1$, the following structure equations
hold
\begin{equation}\label{eq: first system}
\begin{split}
\d \theta_1 & =   \phantom{K_2}\,\theta_1\w\theta_2\,,\\
\d \theta_2 & =             K_2\,\theta_2\w\theta_1\,,\\
\d K_2      &= K_{21}\,\theta_1 + 2(K_2{-}1)\,\theta_2\,.
\end{split}
\end{equation}
(The last equation follows from~$K_1=1$ and~\eqref{eq: Frob conds 2}.)
Moreover, the total differential system for~$a$, $b$, and~$p$ 
simplifies to the system
\begin{equation}\label{eq: da db dp on D_1}
\begin{split}
\d a & = \phantom{{}-\bigl(2b+(K_2-2)(p-1)\bigr)}
         \hbox to 0pt{\hss$(2a+p+1)$}\,\theta_1\,,\\
\d b & =  \phantom{{}-\bigl(2b+(K_2-2)(p-1)\bigr)}
           \hbox to 0pt{\hss${}-2b\,(K_2-1)$}\,\theta_1 
         + \phantom{(2a+p+1)}\hbox to 0pt{\hss$(2b-p+1)$}\,\theta_2\,,\\
\d p &=   {}-\bigl(2b+(K_2-2)(p-1)\bigr)\,\theta_1 + (2a+p+1)\,\theta_2\,,
\end{split}
\end{equation}
which, in view of~\eqref{eq: first system}, is Frobenius.

Correspondingly, on the open set~$D_2$, the following structure equations
hold
\begin{equation}\label{eq: first system dual}
\begin{split}
\d \theta_1 & =             K_1\,\theta_1\w\theta_2\,,\\
\d \theta_2 & =   \phantom{K_1}\,\theta_2\w\theta_1\,,\\
\d K_1      &= 2(K_1{-}1)\,\theta_1 + K_{12}\,\theta_2\,.
\end{split}
\end{equation}
(The last equation follows from~$K_2=1$ and~\eqref{eq: Frob conds 2}.)
Moreover, the total differential system for~$a$, $b$, and~$p$ 
simplifies to the system
\begin{equation}\label{eq: da db dp on D_2}
\begin{split}
\d a & = \phantom{{}-(2b-p+1)}
         \hbox to 0pt{\hss$(2a+p+1)$}\,\theta_1
       - \phantom{\bigl(2a-(K_1-2)(p+1)\bigr)}
          \hbox to 0pt{\hss$2a\,(K_1-1)$}\,\theta_2\,, \\
\d b & =  \phantom{{}-(2b-p+1)\,\theta_1} 
         + \phantom{\bigl(2a-(K_1-2)(p+1)\bigr)}
            \hbox to 0pt{\hss$(2b-p+1)$}\,\theta_2\,,\\
\d p &=   {}-(2b-p+1)\,\theta_1 + \bigl(2a-(K_1-2)(p+1)\bigr)\,\theta_2\,,
\end{split}
\end{equation}
which, in view of~\eqref{eq: first system dual}, is Frobenius.

\begin{remark}[Symmetry of cases]\label{rem: symmetry of cases}
Switching the two families of principal curves switches
the two open sets~$D_1$ and~$D_2$ in~$D$.  Thus, in studying the 
local geometry of the operators~$S$ that satisfy one or the other 
of the structure equations~\eqref{eq: first system} 
or~\eqref{eq: first system dual}, it suffices to study one of the
two cases.  

However, the reader should bear in mind that it is possible for~$D_1$
and~$D_2$ to each be a proper, nonempty subset of~$D$, even though~$D$ 
is connected.  See Example~\ref{ex: D1 D2 proper}.  
However, if~$D$ itself is a $\theta$-rectangular, 
then one of the two is the whole of~$D$.
In the general case, on the (nonempty) overlap~$D_1\cap D_2$, both
structure equations hold and the corresponding differential systems
for~$a$, $b$, and~$p$ simplify even further.
\end{remark}

In conclusion, the calculations made so far have established 
the following result:

\begin{theorem}\label{thm: two types}
If~$S$ is an invertible, nondegenerate operator on a connected
domain~$D$ for which the differential system for~$a$, $b$, and~$p$
is Frobenius, then either~$D$ is the union of two open sets~$D_1$
{\upshape(}on which the equations~\eqref{eq: first system} 
hold\/{\upshape)} and~$D_2$ {\upshape(}on 
which the equations~\eqref{eq: first system dual} hold\/{\upshape)},
or else the equations~\eqref{eq: second system} hold throughout~$D$.

Conversely, if~$S$ is an invertible, nondegenerate operator on
a domain~$D$ for which one of these three systems of structure
equations hold, then the differential system for~$a$, $b$, and~$p$
is Frobenius.
\end{theorem}

\begin{corollary}\label{cor: local existence}
If $S$ is an invertible, nondegenerate operator on a domain~$D$
whose associated coframing~$\theta$ satisfies 
either~\eqref{eq: first system},~\eqref{eq: first system dual},
or~\eqref{eq: second system}, then every point of the domain~$D$
has an open neighborhood on which~$S$ can be realized by a $3$-parameter
family of mutually noncongruent immersions.
\end{corollary}

\begin{proof}
Since the differential system for~$a$, $b$, and~$p$ is Frobenius in any
of these cases, one can choose initial values~$(a_0,b_0,p_0)$ 
of~$a$, $b$, and~$p$ arbitrarily at the specified point 
and have a (unique) solution to the system,
globally defined on a simply connected neighborhood.  
Choosing~$(a_0,b_0,p_0)$ in a bounded region in the quarter-space
defined by~$a_0>0$ and~$b_0>0$ will then yield a congruence class of
immersions realizing~$S$ on a (possibly smaller) fixed neighborhood
of the given point in~$D$.
\end{proof}

\begin{remark}[Locality]\label{rem: locality}
As will be seen in the examples below, the restriction to an open 
neighborhood in Corollary~\ref{cor: local existence} is necessary.  
The point is that,
even if the chosen initial values~$(a_0,b_0,p_0)$ satisfy~$a_0>0$ and~$b_0>0$,
there is no guarantee that the corresponding solutions~$a$ and~$b$
(which, since the $abp$ system is linear, are globally defined if~$D$ 
is simply connected) will be positive throughout~$D$.
\end{remark}

\begin{definition}[The two types]\label{def: two types}
A nondegenerate, invertible operator~$S$ on a domain~$D$ will be
said to be of \emph{Type I} if it satisfies either~\eqref{eq: first system}
or~\eqref{eq: first system dual} and will be said to be of~\emph{Type II}
if it satisfies~\eqref{eq: second system}.
\end{definition}

\section[Integrating the Structure Equations]
{Integrating the Structure Equations}\label{sec: str eqs}

In this final section, the problem of actually integrating the
equations derived in the previous section will be addressed.

\subsection{Operators of Type~I}\label{ssec: Type I}
In this subsection, the operators~$S$ of Type~I will
be studied and it will be shown how to integrate the
equations that define them.  

The definitions of the previous section attach a 
coframing~$\theta=(\theta_1,\theta_2)$ on~$D$ to any 
invertible, nondegenerate operator~$S:TD\to TD$.
The Type~I conditions on~$S$ are then expressed
in terms of this coframing. 

For simplicity, only the cases where either~\eqref{eq: first system}
or~\eqref{eq: first system dual} holds throughout~$D$ will
be considered.  Since these two sub-types differ only in
which family of principal curves is designated first or second,
it suffices to consider only one case.  Thus, it will further 
be assumed that~$S$ satisfies~\eqref{eq: first system}.

\subsubsection{Natural principal coordinates}
\label{sssec: nat prin coords I}
The first task is to find a normal form 
for the coframings~$\theta=(\theta_1,\theta_2)$ on a domain~$D$
that satisfy~\eqref{eq: first system}.  

It is convenient to extend some terminology
to ($2$-dimensional) domains~$D$ endowed with 
a coframing~$\theta$.  The (connected) integral curves of~$\theta_2=0$
will be said to be \emph{principal curves of the first family} 
while the (connected) integral curves of~$\theta_1=0$ will be said to be
\emph{principal curves of the second family}. 
A subdomain~$R\subset D$ will be said to be
\emph{$\theta$-rectangular} if each principal
curve of the first family meets each principal curve
of the second family in a unique point and each principal
curve of the second family meets each principal curve
of the first family in a unique point.  Note that each
point of~$D$ has a rectangular open neighborhood, even
though this may not be computable by quadrature (say).   

\begin{lemma}\label{lem: normal form I}
Suppose that~$\theta=(\theta_1,\theta_2)$ is a coframing
on a domain~$D$ for which there exists a function~$K_2$
so that~$(\theta_1,\theta_2,K_2)$ satisfies~\eqref{eq: first system}.
If~$D$ is~$\theta$-rectangular, 
then there exist functions~$x$ and~$z>0$ on~$D$ so that
\begin{equation}\label{eq: normal form I}
\theta_1 = \frac{\d x}{z}\,,
\qquad\qquad
\theta_2 = \frac{\d z - \bigl(\mu(x)\,z^2+1\bigr)\,\d x}{z}\,,
\qquad\qquad
K_2 = 1 - \mu(x)\,z^2\,,
\end{equation}
where~$\mu$ is a function of a single variable defined on the range of~$x$.  
The functions~$x$ and~$z$ are unique up to a replacement
$(x,z)\mapsto(\lambda\,x+\tau,\,\lambda\,z)$,
where~$\lambda>0$ and~$\tau$ are constants.  

Conversely, if~$\mu$ is a differentiable function defined 
on an interval~$I\subset\bbR$, then the
formulae~\eqref{eq: normal form I} 
define two $1$-forms~$\theta_1$ and~$\theta_2$ and a function~$K_2$ 
on the domain~$D = I\times\R+\subset\R2$ 
that satisfy~\eqref{eq: first system}.
\end{lemma}

\begin{remark}[A quadratic form]\label{rem: a quadratic form}
Note that the quadratic form~$\mu(x)\,(\d x)^2 = (1-K_2)\,{\theta_1}^2$ 
is independent of the choice of local coordinates~$(x,z)$. 
\end{remark}

\begin{proof}
Choose a principal curve from the first family and use its intersections
with the principal curves from the second family as
initial points to construct%
\footnote{This is what is classically known as a `quadrature with
parameters', since it is done by integrating with respect to one
variable while carrying the other along as a `parameter'.  Note
that it is \emph{not} necessary to know the integral curves of~$\theta_1$
explicitly to carry out this construction.} via integration
a function~$w$ on~$D$ so that~$\d w \equiv \theta_2\mod\theta_1$.
Now set~$z = \exp(w)>0$, so that~$(\d z)/z\equiv \theta_2\mod\theta_1$.

Since~$\d \theta_1 = \theta_1\w\theta_2$ (the
first equation of~\eqref{eq: first system}), it follows that
\begin{equation}
\d(z\,\theta_1) = \d z \w \theta_1 + z\,\theta_1\w\theta_2 = 
(\d z - z\,\theta_2)\w \theta_1 = 0.
\end{equation} 
Thus, one can find (by quadrature) a function~$x$ on~$D$
so that~$z\,\theta_1 = \d x$, i.e., so that~$\theta_1 = (\d x)/z$.
The function~$x:D\to\bbR$ is a principal coordinate and
its fibers are the principal curves of the second family.
In particular, the mapping~$(x,z):D\to\R2$ embeds~$D$ as a domain
in~$\R2$ that lies inside~$x(D)\times\R{+}$.

Note that if $(\d x)/z = (\d X)/Z$ for some functions~$X$ and~$Z>0$, 
then~$X = f(x)$ for some function~$f$ of one variable with positive 
derivative and~$Z = f'(x)z$.  Conversely, for any function~$f:x(D)\to\bbR$
with~$f'>0$, the functions~$(X,Z)=\bigl(f(x),f'(x),z\bigr)$ satisfy
the conditions~$(\d Z)/Z\equiv\theta_2\mod\theta_1$ 
and~$\theta_1=(\d X)/Z$.

Now, since~$\d(K_2)\equiv 2(K_2-1)\,\theta_2\mod\theta_1$ (by
the third equation of~\eqref{eq: first system}), it follows that
\begin{equation}
\d \bigl((K_2{-}1)/z^2\bigr) 
\equiv \bigl(2(K_2-1)\,\theta_2\bigr)/z^2 
           + (K_2-1)\bigl((-2/z^2)\,\theta_2\bigr)
\equiv 0 \mod\theta_1\,,
\end{equation}
so there exists a function~$\mu$ defined on~$x(D)\subset\bbR$ so
that~$K_2 = 1 - \mu(x)\,z^2$.  

Next, by construction, $\theta_2-(\d z)/z$ is 
a multiple of~$\theta_1= (\d x)/z$
and, moreover, since~$\d\theta_2 = K_2\,\theta_2\w\theta_1$
(the second equation of~\eqref{eq: first system}),
it follows that the $1$-form $\theta_2-(\d z)/z - (K_2{-}1)\,\theta_1$ 
is closed. Since this $1$-form is closed and a multiple of~$\d x$, 
it must be of the form~$-\nu(x)\,\d x$ 
for some function~$\nu$ defined on~$x(D)\subset\bbR$. Thus,
\begin{equation}
\theta_2 = \frac{\d z - \bigl(\mu(x)\,z^2 + \nu(x)\,z +1\bigr)\,\d x}{z}\,.
\end{equation}

Now, consider any coordinate system~$(X,Z):D\to\R2$ in which
\begin{equation}
\theta_1 = \frac{\d X}{Z}\,,\qquad\qquad\qquad
\theta_2 = \frac{\d Z -
             \bigl(\bar\mu(X)\,Z^2+ \bar\nu(X)\,Z+1\bigr)\,\d X}{Z}
\end{equation}
for some functions~$\bar\mu$ and~$\bar\nu$ on the range of~$X$.
As has been noted already, there is an~$f:x(D)\to X(D)$
so that~$X = f(x)$ and~$Z = f'(x)z$.  Using this
to compare the two formulae for~$\theta_2$ yields
\begin{equation}
\bar\mu\bigl(f(x)\bigr)
 =\frac{\mu(x)}{f'(x)^2}\qquad\text{and}\qquad
\bar\nu\bigl(f(x)\bigr)
 =\frac1{f'(x)}\left(\nu(x)-\frac{f''(x)}{f'(x)}\right)\,.
\end{equation}
Thus, choosing~$f:x(D)\to\bbR$ so that~$f'>0$ 
satisfies~$f''(x) = \nu(x)f'(x)$
(which can be done by a sequence of two quadratures) yields a 
coordinate change~$(X,Z) = \bigl(f(x),f'(x)z\bigr)$ 
for which~$\bar\nu(X)\equiv0$. This yields the normal form of the lemma.

Note that if~$\bar\nu(X)\equiv0$ and~$\nu(x)\equiv0$, 
then~$f''(x)\equiv0$, so that~$f$ is of the 
form~$f(x) = \lambda\,x + \tau$ for some constants~$\lambda>0$ and~$\tau$, 
as claimed in the lemma.

Finally, that the $1$-forms and function defined in~\eqref{eq: normal form I} 
do satisfy~\eqref{eq: first system} can be safely left to the reader.
\end{proof}

\begin{remark}[Computability 2]\label{rem: computability 2}
Note that the desired normal form can be computed by (parametrized)
quadrature.  However, even this step can be eliminated if either~$K_2-1$
is nowhere vanishing or is everywhere vanishing.  

In the first place, if~$K_2-1$ is nowhere vanishing, then the
function~$z>0$ that needs to be found first can simply be taken to be
$z = |K_2-1|^{1/2}$, as this function satisfies the requirement
that~$\d z \equiv z\,\theta_2\mod\theta_1$.  Then a single
quadrature constructs~$x$ so that~$\theta_1 = (\d x)/z$.
This alternative construction provides
coordinates~$(x,z)$ that satisfy~$\theta_1 = (\d x)/z$ 
and~$\theta_2 = \bigl(\d z - (\pm\,z^2{+}\nu(x)\,z{+}1)\,\d x\bigr)/z$, 
instead of the normal form of Lemma~\ref{lem: normal form I}.  The
chief drawback of this normal form is that it cannot be constructed
on a neighborhood of a point where~$K_2{-}1$ vanishes.

On the other hand, if~$K_2-1$ vanishes identically, then,
by the structure equations,~$\theta_1+\theta_2$ is closed and
hence can be written in the form~$(\d z)/z$ by ordinary quadrature, 
thus directly furnishing the desired~$z$.
\end{remark}

\begin{remark}[Type I generality]\label{rem: Type I generality}
One can interpret Lemma~\ref{lem: normal form I} as saying that,
up to local equivalence, the coframings of Type~I depend on one 
arbitrary function of one variable.  It is tempting to regard~$\mu$
as the arbitrary function that `parametrizes' such coframings, but
one must bear in mind that it is not~$\mu$ itself, but the quadratic
form~$\mu(x)\,(\d x)^2$ coupled with the `flat' affine structure
on the space of second principal curves provided 
by Lemma~\ref{lem: normal form I} that provide the distinguishing
invariants.
\end{remark}

Now, the reader will have noticed that~$x$ is a principal coordinate,
but~$z$ is not.  In fact, a second principal coordinate \emph{cannot} 
be constructed by quadrature in general.  However, the following
procedure will `construct' such a coordinate:

Let~$\phi_0$ and~$\phi_1$ be linearly independent solutions on
the interval~$x(D)$ of the second order ordinary differential equation%
\footnote{Note that, for general~$\mu$, these two solutions cannot 
be constructed by quadrature.}
\begin{equation}\label{eq: 2nd order ODE}
\phi''(x) + \mu(x)\,\phi(x) = 0.
\end{equation}
Of course, the Wronskian~$\phi_0\phi'_1-\phi_1\phi'_0$ is constant.
If  $\phi_0\phi'_1-\phi_1\phi'_0 = 1$, the pair~$(\phi_0,\phi_1)$ 
will be said to be \emph{normalized}.  Any two normalized pairs
differ by a unimodular change of basis, and henceforth~$(\phi_0,\phi_1)$
will denote a normalized pair of solutions to~\eqref{eq: 2nd order ODE}
unless it is explicitly stated otherwise.

I claim that there is a solution~$\phi_1$ of~\eqref{eq: 2nd order ODE}
that is positive and increasing on all of~$x(D)$.  To see this, consider
a fixed principal curve~$\Gamma$ in~$D$ of the first family.  Because~$D$
is $\theta$-rectangular,~$\Gamma$ will be mapped by~$(x,z)$ 
to a graph over~$x(D)$ of the form~$z = f(x)$.  Since~$\theta_2$ vanishes
on~$\Gamma$, it follows that~$f:x(D)\to\R+$ must satisfy the equation
\begin{equation}
f'(x) = 1 + \mu(x)\,f(x)^2.
\end{equation} 
Now let~$\phi_1$ be a solution of the linear ODE
\begin{equation}
\phi'_1(x) = \frac1{f(x)}\,\phi_1(x)
\end{equation}
that is positive somewhere (and hence everywhere) on~$x(D)$.  
By construction,
\begin{equation}
\phi''_1(x) = -\frac{f'(x)}{f(x)^2}\,\phi_1(x) 
                + \frac1{f(x)}\,\left(\frac1{f(x)}\,\phi_1(x)\right)
=-\mu(x)\,\phi_1(x)\,,
\end{equation}
and~$\phi_1$ has the desired properties.  Let~$\phi_0$ then be
chosen so that the pair~$(\phi_0,\phi_1)$ is normalized and so
that~$\phi'_0$ vanishes somewhere in~$x(D)$.  Then 
the~$\theta$-rectangularity of~$D$ implies that~$\phi_0(x)-z\,\phi'_0(x)$
is nonvanishing on~$D$.  

Now, define a new function~$y$ on~$D$
by the formula
\begin{equation}\label{eq: y explicit}
y = \frac{\bigl(\phi_1(x) - z\,\phi'_1(x)\bigr)}
         {\bigl(\phi_0(x) - z\,\phi'_0(x)\bigr)}\,.
\end{equation}
This relation can be solved for~$z$ in the form
\begin{equation}\label{eq: y implicit}
z = \frac{\bigl(\phi_1 (x) - y\,\phi_0 (x)\bigr)}
         {\bigl(\phi'_1(x) - y\,\phi'_0(x)\bigr)}\,.
\end{equation}

One then finds that
\begin{equation}\label{eq: 2nd principal coord}
\theta_2 = \frac{ {} - \d y}
           {\bigl(\phi_1 (x) - y\,\phi_0 (x)\bigr)
            \bigl(\phi'_1(x) - y\,\phi'_0(x)\bigr)}\,.
\end{equation}
Thus,~$y$ is the desired second principal coordinate.  Moreover,
it follows that~$y:D\to\bbR$ is a submersion onto an 
interval~$y(D)\subset\bbR$ and that~$(x,y):D\to x(D)\times y(D)$ 
is a diffeomorphism.  Note that, by construction,~$(x,y)$ carries~$\Gamma$
into the segment~$y=0$.  

Note also that the functions~$\bigl(\phi_1(x)- y\,\phi_0(x)\bigr)$
and~$\bigl(\phi'_1(x) - y\,\phi'_0(x)\bigr)$ are positive on~$D$, 
a fact that will be useful below for further constructions.

Principal coordinates~$(x,y)$ found in this manner will be referred to as 
\emph{natural principal coordinates} for the coframing~$\theta$.

\begin{example}[When $\mu$ is constant]\label{ex: mu const}
The reader will find the study of the cases where~$\mu$ is 
constant to be particularly interesting.  When~$\mu$ is a
constant, the formulae~\eqref{eq: normal form I} define a 
coframing of Type~I on the upper half of the $xz$-plane.
In no case is the entire upper half-plane rectangular with
respect to this coframing.  (The reader may enjoy determining
the maximal $\theta$-rectangular subdomains in each case.)
\end{example}

\begin{example}[A locally Type I coframing]\label{ex: D1 D2 proper}
Natural principal coordinates can be used to construct an example of
a connected domain~$D\subset\R2$ with a coframing that is locally
of Type~I, but so that~$D_1$ and~$D_2$ (as defined in 
\S\ref{sssec: special assumption}) 
are each nonempty proper subsets of~$D$.  Of course, 
such a domain cannot be $\theta$-rectangular.  Here is how this 
can be done:

First, let~$\mu:\bbR\to\bbR$ be a smooth function that satisfies~$\mu(x)=0$
for~$x\ge0$ but~$\mu(x)<0$ for~$x<0$.  Let~$\phi_0$ be the function on~$\bbR$
that satisfies~$\phi_0'' + \mu\,\phi_0=0$ and $\phi_0(x) = 1$ for~$x\ge0$
and let~$\phi_1$ be the function on~$\bbR$ that 
satisfies~$\phi_0'' + \mu\,\phi_0=0$ and~$\phi_1(x) = x$ for~$x\ge0$.
Of course, this is a normalized pair~$(\phi_0,\phi_1)$ for the equation
$\phi'' + \mu\,\phi=0$.  Note that because of the sign of~$\mu$ on
the negative reals, $\phi_0$ is decreasing and concave up on the negative
reals while~$\phi_1$ is increasing and concave down on the negative reals.

Now consider the coframing of Type~I
\begin{equation}\label{eq: half-flat example}
\theta_1 = \frac{\bigl(\phi'_1(x) - y\,\phi'_0(x)\bigr)\,\d x}
                {\bigl(\phi_1 (x) - y\,\phi_0 (x)\bigr)}\,,\qquad\qquad
\theta_2 = \frac{ {} - \d y}
           {\bigl(\phi_1 (x) - y\,\phi_0 (x)\bigr)
            \bigl(\phi'_1(x) - y\,\phi'_0(x)\bigr)}\,,
\end{equation}
that is smooth and well-defined on the open domain~$D_1\subset\R2$
defined by the inequalities~$y<0$ when~$x\ge0$ and
\begin{equation}
\frac{\phi'_1 (x)}{\phi'_0 (x)} < y < \frac{\phi_1 (x)}{\phi_0 (x)}
\end{equation}
when~$x<0$.  Note that, when~$y<0<x$, i.e., in the fourth quadrant
of the plane, the above formulae simplify to
\begin{equation}
\theta_1 = \frac{\d x}{x-y}\,,\qquad\qquad
\theta_2 = \frac{\d y}{y-x}\,.
\end{equation}

Thus, the involution~$\Phi:\R2\to\R2$ defined by~$\Phi(x,y) = (-y,-x)$,
which preserves the fourth quadrant, satisfies~$\Phi^*\theta_1=\theta_2$
and~$\Phi^*\theta_2=\theta_1$ there. 

Finally, let~$D_2 = \Phi(D_1)$ and extend~$\theta_1$ and~$\theta_2$ 
to~$D = D_1\cup D_2 = \Phi(D)$ in the obvious way so that~$\Phi^*$
exchanges~$\theta_1$ and~$\theta_2$ globally on~$D$.  

Since~$\mu$ is nonzero on the negative real axis, the function $K_1-1$
is nonvanishing when~$y>0$ and the function~$K_2-1$ is nonvanishing
when~$x<0$.  Thus, this is the desired example.
\end{example}

\subsubsection{Integrals of the Frobenius system}\label{sssec: Frob ints I}
The integrals of the system for~$a$, $b$, and~$p$ are easily described
in the coordinates~$(x,z)$ of Lemma~\ref{lem: normal form I}.  One
finds that the following formulae hold
\begin{equation}\label{eq: abp integrated I}
\begin{split}
a &= -1+f(x)\\
b &= -\bigl(1+\mu(x)z^2\bigr)f(x)+z\,f'(x)-{\ts\frac12}z^2\,f''(x)\\
p &= 1-2\,f(x)+z\,f'(x)
\end{split}
\end{equation}
where~$f$ is any solution of the equation
\begin{equation}\label{eq: 3rd for f}
f'''(x) + 4\,\mu(x)\,f'(x) + 2\mu'(x)\,f(x) = 0.
\end{equation}
The general solution of this equation is easily seen to be
\begin{equation}
f(x) = {} - c_0\,\phi_0(x)^2 - 2c_1\,\phi_0(x)\phi_1(x) - c_2\,\phi_1(x)^2
\end{equation}
where~$(\phi_0,\phi_1)$ is a normalized pair of solutions 
of~\eqref{eq: 2nd order ODE} and~$c_0$, $c_1$, and~$c_2$ are arbitrary 
constants.

In terms of natural principal coordinates as described above,
the formulae for~$a$ and~$b$ simplify%
\footnote{The corresponding formula for~$p$ is not as simple, 
but will not be needed in what follows.} 
to
\begin{equation}\label{eq: ab in nat prin coords}
\begin{split}
a &= -1-c_0\,\phi_0(x)^2-2c_1\,\phi_0(x)\phi_1(x)-c_2\,\phi_1(x)^2\,,\\
b &=  \frac{(c_0 + 2c_1\,y + c_2\,y^2)}
        {\bigl(\phi'_1(x)-y\,\phi'_0(x)\bigr)^2}\,.
\end{split}
\end{equation}
In this form, it is not difficult to understand how to choose the
constants~$c_i$ so that~$a$ and~$b$ will be positive at a given point
of~$D$.  In particular, these constants must satisfy~${c_1}^2-c_0\,c_2>0$
or else it will be impossible for~$a$ and~$b$ to be positive simultaneously.

Note also that, because~$z$ is strictly positive on~$D$, the
expression~$\phi_1(x)-y\,\phi_0(x)$ cannot vanish.  This implies that,
at any point of~$D$, the allowable values of~$(c_0,c_1,c_2)$ for which
$a$ and~$b$ will be positive at the specified point consists of the
(non-empty) intersection of two open half-spaces (with non-parallel
bounding planes).  

Conversely, it is not difficult to see that,
for any given values of~$c_0$, $c_1$, and~$c_2$, the set of points of~$D$
at which both~$a$ and~$b$ are positive is a disjoint union of open
rectangles in~$D$.  (Of course, this uses the assumption that~$D$ itself
is $\theta$-rectangular.)

\subsubsection{Recovering~$S$}\label{sssec: recovering S}
So far, the discussion in this section has shown how one can
write down $1$-forms~$\theta_1$ and~$\theta_2$ 
and a function~$K_2$ on a domain~$D$ satisfying~\eqref{eq: first system}.  
However, it is not immediate whether or not such a system necessarily 
comes from a nondegenerate, invertible operator~$S$ defined on~$D$,
and, if so, `how many' such operators~$S$ there are.  It is now 
time to address this question.

By definition, the desired~$S$, if it exists, will be of the form
\begin{equation}\label{eq: S in UVtheta}
S = \frac1U\,\,\tb_1\otimes\theta_1 + \frac1V\,\,\tb_2\otimes\theta_2\,,
\end{equation}
where~$\tb_1$ and~$\tb_2$ are the vector fields on~$D$ dual to the
coframing defined by the $1$-forms~$\theta_1$ and~$\theta_2$ and
where~$U$ and~$V$ are nonzero and nowhere equal functions on~$D$ 
that satisfy
\begin{equation}\label{eq: UV pde in theta}
\begin{split}
\d U &= U_1\,\theta_1 + (U-V)\,\theta_2\,,\\
\d V &= (V-U)\,\theta_1 + V_2\,\theta_2\,,
\end{split}
\end{equation}
for some functions~$U_1$ and~$V_2$ on~$D$.  Conversely,
by the very definitions of~$\theta_1$ and~$\theta_2$, 
if~$U$ and~$V$ are nonvanishing, nowhere equal functions on~$D$ that 
satisfy~\eqref{eq: UV pde in theta} for some functions~$U_1$ and~$V_2$,
then~\eqref{eq: S in UVtheta} defines an invertible, nondegenerate
operator on~$TD$ that is of Type~I.

The system~\eqref{eq: UV pde in theta} constitutes a pair of linear,
first order partial differential equations for~$(U,V)$.  In fact, this
is a hyperbolic system whose characteristics are the principal curves,
i.e., the level curves of~$x$ and~$y$.  For example, if~$(x,y):R\to\R2$ 
are principal coordinates on a $\theta$-rectangle~$R\subset D$, then
there exist nonvanishing functions~$s$ and~$t$ on~$(x,y)(R)$ so that
\begin{equation}\label{eq: theta in xy general}
\theta_1 = s(x,y)\,\d x\,,\qquad\qquad\qquad
\theta_2 = t(x,y)\,\d y\,.
\end{equation}
In these local coordinates, the equations~\eqref{eq: UV pde in theta}
become the coupled linear system of PDE
\begin{equation}\label{eq: UV pde in xy}
\frac{\partial U}{\partial y} = t(x,y)\,(U-V)\,,\qquad\qquad\qquad
\frac{\partial V}{\partial x} = s(x,y)\,(V-U)\,.
\end{equation}
This system is visibly hyperbolic, with~$x$ and~$y$ being the
characteristic directions.   Standard existence theorems
ensure that there exist solutions, locally.  In fact, one can
specify~$U$ and~$V$ arbitrarily along a noncharacteristic curve~$\Gamma$ 
in~$R$ (i.e., a curve that is everywhere transverse to the principal 
curves) and there will be an neighborhood of~$\Gamma$ in~$R$ on which
a solution to~\eqref{eq: UV pde in xy} exists and assumes the 
prescribed values on~$\Gamma$.  In particular, if one specifies~$U$
and~$V$ so that they are unequal and nonvanishing along~$\Gamma$,
then the pair~$(U,V)$ will be a solution with the desired properties.%
\footnote{The reader will note that this discussion of recovering
$U$ and~$V$ from $\theta$ does not depend on $\theta$ satisfying
the conditions for Type~I (or any conditions, for that matter).
Thus, any coframing~$\theta$ is locally realizable as the coframing
of some invertible, nondegenerate operator~$S$.}

However, it is not necessary to appeal to such theorems to prove 
existence.  It turns out that the system~\eqref{eq: UV pde in theta}
is integrable by the method of Darboux, as will now be explained.

Assume that a solution~$(U,V)$ to~\eqref{eq: UV pde in theta} exists
and compute the exterior derivatives of the 
equations~\eqref{eq: UV pde in theta}, using the structure 
equations~\eqref{eq: first system} and the 
equations~\eqref{eq: UV pde in theta} themselves.  
The result can be written in the form

\begin{equation}\label{eq: U1V2 relations}
\begin{split}
0 &= \left(\d U_1 - 2U_1\,\theta_2 + (K_2-1)\,\d U \right)\w\theta_1\,,\\
0 &= \left(\d V_2 -(K_2+1)V_2\,\theta_1\right)\w\theta_2\,,
\end{split}
\end{equation}

The first of these equations is just
\begin{equation}
\d U_1 -2U_1\, \theta_2  + (K_2-1)\,\d U \equiv 0 \mod \theta_1
\end{equation}

Using the coordinates~$(x,z)$ guaranteed 
by Lemma~\ref{lem: normal form I}, this equation takes the form
\begin{equation}
\d U_1 -2\,U_1\, \frac{\d z}{z}  -\mu(x)\,z^2\,\d U \equiv 0 \mod \d x.
\end{equation}
Dividing this equation by~$z^2$, it can be rewritten in the form
\begin{equation}
\d \left(\frac{U_1}{z^2}-\mu(x)\, U\right) \equiv 0 \mod \d x.
\end{equation}
Thus, there must exist a function~$\xi$ defined on~$x(D)$ 
such that
\begin{equation}\label{eq: U1 as xi and U}
U_1 = \bigl(\mu(x)\,U + \xi(x)\bigr)\,z^2.
\end{equation}

For the second equation of~\eqref{eq: U1V2 relations}, if one
uses the local normal form  to expand the right hand side
and makes the substitutions
\begin{equation}\label{eq: V2 z substitutions}
z = \frac{\bigl(\phi_1 (x) - y\,\phi_0 (x)\bigr)}
         {\bigl(\phi'_1(x) - y\,\phi'_0(x)\bigr)}\,,
\qquad\qquad
V_2 = Q\,\bigl(\phi_1(x)-y\,\phi_0(x)\bigr)^2\,
      \bigl(\phi'_1(x)-y\,\phi'_0(x)\bigr)\,,
\end{equation}
for some new variable~$Q$, 
then this equation becomes
\begin{equation}
0 = - \bigl(\phi_1(x)-y\,\phi_0(x)\bigr)\,(\d Q \w \d y).
\end{equation}
Since the scalar factor in this equation is nonvanishing on~$D$,
it follows that~$\d Q\w\d y=0$.  In other words, there is a 
function~$\eta$ defined on~$y(D)$ such that
\begin{equation}\label{eq: V2 as eta}
V_2 = \bigl(\phi_1(x)-y\,\phi_0(x)\bigr)^2\,
      \bigl(\phi'_1(x)-y\,\phi'_0(x)\bigr)\,\eta(y).
\end{equation}

Substituting the relations~\eqref{eq: U1 as xi and U} 
and~\eqref{eq: V2 as eta} back into~\eqref{eq: UV pde in theta}
then yields a differential system for~$U$ and~$V$ that is
Frobenius for any choice of~$\xi$ and~$\eta$.  In fact, this
system can be integrated explicitly in the form
\begin{equation}\label{eq: UV explicit}
\begin{split}
U & = f(x) - \frac{g(y)}{\phi'_1(x)-y\,\phi'_0(x)}
     - \left(\frac{{\phi_1(x)-y\,\phi_0(x)}}
                   {{\phi'_1(x)-y\,\phi'_0(x)}}\right)\,f'(x) \,,\\[10pt]
V & = f(x) - \phi_0(x)\,g(y) - \bigl(\phi_1(x)-y\,\phi_0(x)\bigr)\,g'(y)\,,
\end{split}
\end{equation}
where~$f$ and~$g$ satisfy the equations%
\footnote{
Of course,~$f$ and~$g$ can be computed from~$\xi$ and~$\eta$ by quadratures
since~$\phi_0$ and~$\phi_1$ are assumed known.
} 
\begin{equation}\label{eq: g h quadratures}
f''(x) + \mu(x)\,f(x) = -\xi(x)\,,\qquad\qquad g''(y) = \eta(y)\,.
\end{equation}

Conversely, for any smooth functions~$f$ on~$x(D)$ and~$g$ on~$y(D)$,
the formulae~\eqref{eq: UV explicit} define a solution 
of~\eqref{eq: UV pde in theta}.  Thus, it follows that~\eqref{eq: UV explicit}
is the general solution of~\eqref{eq: UV pde in theta}.

\begin{theorem}\label{thm: Type I explicit solution}
If~$S:TD\to TD$ is an invertible, nondegenerate operator 
whose~$\theta$-coframing
satisfies the Type~I structure equations~\eqref{eq: first system}, then 
each point of~$D$ has a $\theta$-rectangular neighborhood~$R$ on which 
there exist principal coordinates~$(x,y)$ in which~$S$ takes the form
\begin{equation}\label{eq: explicit S}
S =  \frac1U\,\,\frac{\partial\hfil}{\partial x}\otimes \d x
   + \frac1V\,\,\frac{\partial\hfil}{\partial y}\otimes \d y
\end{equation}
where~$U$ and~$V$ are of the form~\eqref{eq: UV explicit} for some
functions~$g$ on~$y(R)$ and~$f$, $\phi_0$, and~$\phi_1$ on~$x(R)$ 
satisfying~$\phi_0\phi'_1-\phi_1\phi'_0=1$.  

Moreover, the coordinates~$x$ and $y$ and the functions~$g$, $f$, $\phi_0$
and~$\phi_1$ are computable from~$S$ by algebraic operations, quadratures,
and the integration of a single, linear, self-adjoint second order ordinary
differential equation on~$x(D)$. 

Conversely, if~$R$ is a rectangle in the $xy$-plane, then, for any choice 
of $g$ on~$y(R)$ and~$f$, $\phi_0$, and~$\phi_1$ 
on~$x(R)$ satisfying~$\phi_0\phi'_1-\phi_1\phi'_0=1$ such that
that the functions~$U$ and $V$ defined by~\eqref{eq: UV explicit} are
nonvanishing and unequal on~$R$, the formula~\eqref{eq: explicit S}
defines an invertible, nondegenerate operator~$S:TR\to TR$ 
whose $\theta$-coframing satisfies~\eqref{eq: first system} and for
which~$(x,y)$ is a $\theta$-principal coordinate system. \qed
\end{theorem}

\begin{remark}[Normal form ambiguities]\label{rem: normal form ambiguities}
The coordinates~$x$ and $y$ and the functions~$g$, $f$, $\phi_0$, 
and~$\phi_1$ are not quite canonically determined by~$S$: 

The coordinate~$x$ is determined up to an affine transformation~$x\mapsto 
\lambda\,x{+}\tau$ where~$\lambda>0$ and~$\tau$ are constants. 

Once~$x$ is chosen, the function~$\mu$ can be found and
then the normalized pair $(\phi_0,\phi_1)$ is determined up to a
(constant) unimodular change of basis.  

Once~$x$, $\phi_0$ and~$\phi_1$ are chosen, the function~$y$ is
determined.  (It may be necessary to re-choose the normalized 
pair~$(\phi_0,\phi_1)$ so that~$y$ remains finite on all of~$R$.)

Finally, the functions~$f$ and~$g$ are determined up to a 
replacement of the form
\begin{equation}
\bigl(f(x),g(y)\bigr) 
= \bigl(f(x){+}c_1\,\phi_1(x){+}c_0\,\phi_0(x),\,g(y){+}c_1\,y{+}c_0\,\bigr)
\end{equation}   
for any two constants~$c_0$ and~$c_1$.
\end{remark}

\begin{remark}[Reduction to quadrature]\label{rem: reduction to quadrature}
As Theorem~\ref{thm: Type I explicit solution} shows, the entire process 
of computing a normal form for~$S$ requires only algebraic operations, 
quadratures, and the solution of a single second, order self-adjoint
ordinary differential equation, namely~\eqref{eq: 2nd order ODE}.  

However, given~$S$, one can usually dispense with this last step, 
since an alternative is available.  In particular, as long as~$V_2$ 
is not identically vanishing (which is the same as 
the condition~$\d''V\not\equiv0$), this can be done as follows: 

In terms of principal coordinates, 
the formulae~\eqref{eq: half-flat example}, 
\eqref{eq: V2 z substitutions}, and~\eqref{eq: V2 as eta}
show that
\begin{equation}\label{eq: V2 quadrature}
\left(\d''V\right)^2\circ\theta_2\circ\theta_1 
 = {V_2}^2\,{\theta_2}^3\circ\theta_1 
 = {} - \eta(y)^2\,(\d y)^3\circ \d x\,.
\end{equation} 
Of course, the quartic form on the left hand side is computable from~$S$ 
by differentiation alone.  

Assume that the left hand side of~\eqref{eq: V2 quadrature}
is nonzero.  Then by algebraic operations, one can write it 
in the form~$\psi_1\circ{\psi_2}^3$ for some coframing~$(\psi_1,\psi_2)$
that is unique up to a replacement of the 
form
\begin{equation}\label{eq: psi ambig}
(\psi_1,\psi_2)\mapsto (r^{-3}\,\psi_1,r\,\psi_2)
\end{equation}
for some nonvanishing function~$r$ on~$D$.  By differentiation, one can
now find a unique $1$-form~$\rho$ satisfying~$\d\psi_1 = -3\,\rho\w\psi_1$
and~$\d\psi_2=\rho\w\psi_2$ (this~$\rho$ is essentially a connection
form for the quartic).  Under the replacement~\eqref{eq: psi ambig},
one finds that~$\rho$ is replaced by~$\rho+\d(\log r)$.  
Now, by~\eqref{eq: V2 quadrature}, it is obvious that there exists
a choice, namely
\begin{equation}
(\bar\psi_1,\bar\psi_2) 
  = \bigl(\d x,\, -\bigl(\eta(y)\bigr)^{2/3}\,\d y\bigr)
\end{equation} 
for which each~$\bar\psi_i$ is closed, i.e., for which the
corresponding connection form is~$\bar\rho=0$.  
Consequently,~$\rho$ must be closed 
for \emph{any} choice of~$(\psi_1,\psi_2)$.  Since~$\rho$ is 
closed, it can, by quadrature, be written in the form~$\rho = \d(\log r)$ 
for some positive function~$r$ on~$D$.  Then replacing~$(\psi_1,\psi_2)$
by~$(r^3\,\psi_1, r^{-1}\,\psi_2)$ yields a $\psi$-coframing for which
each~$\psi_i$ is closed.  Thus, by quadrature, one
can write~$\psi_1 = \d x$ and~$\psi_2 = \d t$ for some principal
coordinates~$x$ and~$t$ on~$D$.  (I am using~$t$ instead of~$y$ here
because~$t$ will not, in general, be a natural principal coordinate,
although~$x$ certainly is.)  Regarding~$x$ as now known, $z$ can
be defined by requiring that~$\theta_1 = (\d x)/z$.  (It may be necessary
to replace~$x$ by~$-x$ to ensure that~$z$ is positive.)

Now, the curve~$t = t_0$ for some constant~$t_0$
is an integral curve of~$\theta_2$.
In the~$xz$ coordinates, it can be written in the form~$z = f(x)$
for some function~$f$ that necessarily satisfies~$f' = 1 + \mu\,f^2$.  
Of course, as has already been explained, the function~$\phi_1$ that 
satisfies~$\phi'_1 = (1/f)\,\phi_1$ can now be computed by quadrature 
and satisfies~$\phi''_1 + \mu\,\phi_1=0$.  Thus, one has one solution 
of~\eqref{eq: 2nd order ODE} computed by quadrature.  Solving
the first order inhomogeneous equation~$\phi_0\phi'_1-\phi_1\phi'_0 =1$
for~$\phi_0$ by quadrature then yields 
the desired normalized pair~$(\phi_0,\phi_1)$.

Note that this procedure produces \emph{both} principal coordinates
by algebraic operations and simple quadratures.  Its main disadvantage
is that it depends on properties of the pair~$(U,V)$ and not just
on the coframing~$\theta$.
\end{remark}

\begin{remark}[Primary invariants]\label{rem: primary invariants}
The reader will have noticed that there are three fundamental 
invariants that are computable purely by differentiation that,
in some sense, define the principal coordinates up to quadrature:
\begin{equation}\label{eq: three invariants}
\begin{split}
(1-K_2)\,{\theta_1}^2 &= \mu(x)\,(\d x)^2\,,\\
\bigl(U_1-(1{-}K_2)U\,\bigr)\,{\theta_1}^2 &= \xi(x)\,(\d x)^2\,,\\
-{V_2}^2\,{\theta_2}^3\circ\theta_1 
 &= {} \eta(y)^2\,(\d y)^3\circ \d x\,.
\end{split}
\end{equation}
These differential invariants were originally found 
by the method of Darboux, though the treatment above that generates
them was developed so that the reader need not be familiar
with this method.

It can be shown that these invariants plus the affine structure on 
the space of principal curves of the second family (which requires 
a quadrature to compute) are sufficient to test whether or not two 
given operators of Type~I are equivalent up to a change of variable.
\end{remark}

\begin{remark}[A final elimination]\label{rem: final elimination}
Informally, Theorem~\ref{thm: Type I explicit solution} says that
the operators of Type~I depend on three arbitrary functions of 
one variable, i.e., four arbitrary functions~$f(x)$, $\phi_0(x)$, $\phi_1(x)$,
and~$h(y)$ subject to the differential 
equation~$\phi_0\phi'_1{-}\phi_1\phi'_0=1$ and some open conditions.  

It is worth pointing out that one can make a change of variables to
eliminate the differential equation~$\phi_0\phi'_1{-}\phi_1\phi'_0=1$: 
Use the fact that~$\bigl(\phi_0(x),\phi_1(x)\bigr)$ never
vanishes to write
\begin{equation}
\phi_0(x) = r(t)\,\cos t\,,\qquad\qquad\qquad
\phi_1(x) = r(t)\,\sin t
\end{equation}
for some functions~$t$ and~$r>0$ on~$x(R)$.  The 
identity
\begin{equation}
\d x = \bigl(\phi_0(x)\phi'_1(x)-\phi_1(x)\phi'_0(x)\bigr)\,\d x
     = \phi_0\,\d\phi_1-\phi_1\,\d\phi_0 = r^2\,\d t
\end{equation}
shows that~$\d t$ is nonvanishing on~$x(R)$ and hence can be taken
as a coordinate.  (Note that~$t$ is a principal coordinate replacing~$x$.)
In particular, there exists a function~$\rho:t(R)\to\bbR$ 
so that~$r = \rho(t)$.

Now, writing~$g(x) = \gamma(t)$ for some function~$\gamma$ on~$t(R)$,
one computes~$g'(x) = \dot\gamma(t)/\rho(t)^2$.  Similarly,~$\phi'_0(x)$
and~$\phi'_1(x)$ can be expressed in terms of~$t$, $\rho(t)$, 
and~$\dot\rho(t)$.

Thus, all of the expressions in the formulae for~$U$ and~$V$ involving~$x$ 
can be replaced by expressions in~$t$ and the (arbitrary positive) 
function~$\rho$ on~$t(R)$.  

This change of variables expresses~$S$ in $ty$-coordinates 
in terms of the three arbitrary functions~$\rho(t)$, $\gamma(t)$, 
and~$g(y)$ and their derivatives,
where these three functions are only subject to open conditions, not
equations (differential or otherwise).  However, this formula does not
appear to be particularly useful, so it will not be explored further.
\end{remark}

\subsubsection{Integrating the structure equations}
\label{sssec: int str eqs}
Now that explicit formulae have been found for~$\theta_1$,
$\theta_2$, $a$, $b$, $U$ and~$V$, the structure forms~$\omega_1$,
$\omega_2$, $\omega_{31}$, $\omega_{32}$ and~$\omega_{12}$ 
can be regarded as known.  Since, by construction, these forms 
satisfy the structure equations, Bonnet's theorem can be used to
show that there exists a corresponding realization~$\xb$ of
the operator~$S$.  

However, Bonnet's theorem is not an 
effective theorem in the sense that the realization~$\xb$
cannot, in general be computed from the structure forms using
only algebraic operations and quadratures.  The goal of
this subsection is to indicate that, in fact, one can 
avoid having to quote Bonnet's theorem by following an algorithm
for constructing~$\xb$ that only involves algebraic operations
and quadratures.  In the final subsubsections of this section,
this algorithm will be used to compute some explicit examples.

Now, the structure forms~$\omega_{31}$, $\omega_{32}$ and~$\omega_{12}$
have the expressions
\begin{equation}\label{eq: explicit omega_ij}
\begin{split}
\omega_{31} 
   &= \frac{\theta_1}{\sqrt{a}} 
 = \frac{\bigl(\phi'_1(x)-y\,\phi'_0(x)\bigr)\,\d x}
      {\bigl(\phi'_1(x)-y\,\phi'_0(x)\bigr)
      \sqrt{-1-c_0\,\phi_0(x)^2
            -2c_1\,\phi_0(x)\phi_1(x)-c_2\,\phi_1(x)^2\strut}}\,,\\
\omega_{32} 
   &=\frac{\theta_2}{\sqrt{b}}
 =\frac{-\,\d y}
      {\bigl(\phi'_1(x)-y\,\phi'_0(x)\bigr)
         \sqrt{c_0 + 2c_1\,y + c_2\,y^2\strut}}\,,\\
\omega_{12} &= -\frac{\sqrt{b}}{\sqrt{a}}\,\theta_1
               +\frac{\sqrt{a}}{\sqrt{b}}\,\theta_2\\
    &= \frac{-1}{\bigl(\phi'_1(x)-y\,\phi'_0(x)\bigr)}
\left(
\sqrt{\frac{c_0 + 2c_1\,y + c_2\,y^2\strut}
  {-1-c_0\,\phi_0(x)^2
         -2c_1\,\phi_0(x)\phi_1(x)-c_2\,\phi_1(x)^2\strut}}\,\d x
\right.\\
&\qquad\qquad\qquad\qquad\qquad\qquad
+ \left.\sqrt{
\frac{{-1-c_0\,\phi_0(x)^2-2c_1\,\phi_0(x)\phi_1(x)-c_2\,\phi_1(x)^2\strut}}
     {c_0 + 2c_1\,y + c_2\,y^2\strut}}\,\d y
\right).
\end{split}
\end{equation}

These equations already have an interesting consequence:  
Let~$(\eb_1,\eb_2,\eb_3):D\to \Or(3)$ be a solution to the equations
$\d\eb_i = \eb_j\,\omega_{ji}$.  Since~$\omega_{31}\w\omega_{32}\not=0$,
the map~$\eb_3:D\to S^2$ is an immersion. Of course, $(\eb_1,\eb_2)$ is
a tangential orthonormal frame field for this immersion.  Note that,
by construction,~$\sqrt{a}$ is the function that gives the geodesic
curvature of the principal curves from the second family, i.e., the 
integral curves of~$\theta_1=0$.  Now, since~$\d a\equiv 0 \mod \theta_1$,
it follows that the geodesic curvature of the $\eb_3$-image of each
of these curves is constant.  Of course, this implies that the $\eb_3$-image
of each principal curve of the second family is a geodesic circle.
This observation has the following consequence:

\begin{proposition}\label{prop: spherical circles}
Let~$S:TD\to TD$ be an invertible, nondegenerate operator on a domain~$D$
that satisfies~\eqref{eq: first system}.  Then
for any normally oriented realization~$\xb:R\to\E3$ of~$S$
on a subdomain~$R\subset D$, the Gau{\ss} image of each
principal curve of the second family is an arc of a geodesic circle 
on the $2$-sphere.

Moreover, the spherical centers of these geodesic circles lie 
on a fixed great circle on~$S^2$.
\end{proposition}

\begin{proof}
All but the last statement has been verified already.  To prove the
last statement, it suffices to note that the spherical center of the
geodesic circle tangent to~$\eb_2$ with geodesic curvature~$\sqrt{a(x)}$ 
is given by
\begin{equation}\label{eq: z spherical center}
\zb = \frac{-1}{1+a(x)}\,\eb_1 + \frac{\sqrt{a(x)}}{1+a(x)}\,\eb_3\,.
\end{equation}
Computation now shows that $\d\zb = \wb\,\sigma$ where $\wb\cdot\wb=1$ and
\begin{equation}\label{eq: sigma defined}
\sigma = \frac{\sqrt{{c_1}^2-c_0c_2}}{\bigl(1+a(x)\bigr)\sqrt{a(x)}}\,\d x\,,
\end{equation}
and, moreover, that $\d\wb = -\zb\,\sigma$.
It follows that~$\zb$ moves on the great circle perpendicular to
the fixed vector~$\zb\times\wb$.
\end{proof}

Using Proposition~\ref{prop: spherical circles}, it is now not difficult
to integrate the equations and find the~$\eb_i$ explicitly.  In fact,
one does the following:  First, construct, by quadrature, a function~$s(x)$
so that
\begin{equation}\label{eq: s quadrature}
\d s =\sigma=\frac{\sqrt{{c_1}^2-c_0c_2}}
                  {\bigl(1+a(x)\bigr)\sqrt{a(x)}}\,\d x\,.
\end{equation}
One can then show that, up to a rotation,~$\eb_3$ is given by
\begin{equation}\label{eq: integrated e3 I}
\eb_3(x,y) = \left[
\,\begin{matrix}
\ds\cos\bigl(s(x)\bigr)\,\frac{\sqrt{a(x)}}{\sqrt{1+a(x)}}
      -\sin\bigl(s(x)\bigr)\,
\frac{\bigl(c_0\,\phi_0(x)
       +c_1\,\bigl(\phi_1(x)+y\,\phi_0(x)\bigr)+c_2\,y\,\phi_1(x)\bigr)}
{\bigl(\phi_1(x)+y\,\phi_0(x)\bigr)\sqrt{1+a(x)}\sqrt{{c_1}^2-c_0c_2}}\\[15pt]
\ds\frac{\sqrt{c_0 + 2c_1\,y + c_2\,y^2}}
     {\bigl(\phi_1(x)-y\,\phi_0(x)\bigr)\sqrt{{c_1}^2-c_0c_2}}\\[15pt]
\ds\sin\bigl(s(x)\bigr)\,\frac{\sqrt{a(x)}}{\sqrt{1+a(x)}}
      +\cos\bigl(s(x)\bigr)\,
\frac{\bigl(c_0\,\phi_0(x)
       +c_1\,\bigl(\phi_1(x)+y\,\phi_0(x)\bigr)+c_2\,y\,\phi_1(x)\bigr)}
{\bigl(\phi_1(x)+y\,\phi_0(x)\bigr)\sqrt{1+a(x)}\sqrt{{c_1}^2-c_0c_2}}
\end{matrix}\,\right].
\end{equation}
There are, of course, similar formulae for~$\eb_1$ and~$\eb_2$, but
they will not be needed explicitly in the algorithm to be considered.
The important point is that~$\eb_3$ can be constructed by quadrature
and that~$\omega_{31}$ is a multiple of~$\d x$ while $\omega_{32}$
is a multiple of~$\d y$.  This implies that the comparison
\begin{equation}
   \frac{\partial\eb_3}{\partial x}\,\d x 
 + \frac{\partial\eb_3}{\partial y}\,\d y
 = \d \eb_3 = {} -\eb_1\,\omega_{31} -\eb_2\,\omega_{32}
\end{equation}
yields
\begin{equation}
\frac{\partial\eb_3}{\partial x}\,\d x = -\eb_1\,\omega_{31}
\qquad\qquad\text{and}\qquad\qquad
\frac{\partial\eb_3}{\partial y}\,\d y = -\eb_1\,\omega_{32}\,.
\end{equation}

Finally, the immersion~$\xb:D\to\E{3}$ satisfies the structure equation
\begin{equation}\label{eq: dx formula I}
\begin{split}
\d\xb &= \eb_1\,\omega_1 + \eb_2\,\omega_2 
      = U(x,y)\,\eb_1\,\omega_{31} + V(x,y)\,\eb_2\,\omega_{32}\\
      &= {}-U(x,y)\,\frac{\partial\eb_3}{\partial x}\,\d x
           -V(x,y)\,\frac{\partial\eb_3}{\partial y}\,\d y.
\end{split}
\end{equation}
The structure equations imply
that the vector-valued differential form on the right hand
side of~\eqref{eq: dx formula I} is indeed a closed $1$-form, 
so that~$\xb$ can be recovered by quadrature:
\begin{equation}\label{eq: x quadrature general I}
\xb = - \int U(x,y)\,\frac{\partial\eb_3}{\partial x}\,\d x
           +V(x,y)\,\frac{\partial\eb_3}{\partial y}\,\d y\,.
\end{equation}

Thus, the final result is

\begin{theorem}\label{thm: type 1 realization by quadrature}
The local realizations of an operator~$S$ of Type~I can be computed
by quadratures, once the principal coordinates are found.
\end{theorem}

\subsubsection{The case~$\mu=0$}\label{sssec: mu=0}
The rest of this subsection about the geometry of Type~I
operators and their realizations will concern only the
special case~$\mu\equiv0$, i.e., when~$K_2\equiv1$.  This
simplest case has special features that are not shared by 
the general Type~I operators.  For example, some of the
quadratures that are needed in the general case can be
eliminated, thus leading to more explicit formulae and
easily computed examples.   Moreover the domains of the
realizations are more easily described.

Some of the features discussed
here can be generalized to other values of~$\mu(x)$, 
but that will be left to the interested reader.

To distinguish this case from the general case, upper case
letters will be used for the natural principal coordinates.
Thus,~$X$ and~$Z$ instead of~$x$ and~$z$.  

The equation~$\phi''(X)+\mu(X)\,\phi(X)=0$ now simplifies
to~$\phi''(X) = 0$, with the obvious normalized pair~$(\phi_0,\phi_1)
=(1,X)$.  This gives~$Y = X-Z>0$ and the formulae for
the $\theta$-coframing assume the simple, symmetric form
\begin{equation}\label{eq: theta in XY}
\theta_1 = \frac{\d X}{X-Y}\,,\qquad\qquad
\theta_2 = \frac{\d Y}{Y-X}\,,
\end{equation}

The formulae for~$a$, $b$, and~$p$ simplify to
\begin{equation}\label{eq: abp in XY-terms I}
\begin{split}
a &= {\textstyle{\frac12}}(-1 + c_0 + 2c_1\,X + c_2\,X^2)\,,\\
b &= {\textstyle{\frac12}}(-1 - c_0 - 2c_1\,Y - c_2\,Y^2)\,,\\
p &= c_0 + c_1(X+Y) + c_2\,XY\,.
\end{split}
\end{equation}
(In the interests of preserving the~$XY$ symmetry, the
usage of the constants~$c_i$ is now slightly different 
from that of the general case.)

The reader must keep in mind that the conditions~$a>0$
and~$b>0$ must still be imposed.  These inequalities, together with
the requirement~$X>Y$, impose inequalities on~$c_0$,~$c_1$,~and~$c_2$.

These inequalities amount to the condition that there exist constants
~$\xi$, $\eta$, and~$\lambda$ with~$\xi>\eta$ and~$|\lambda|<1$ so that
\begin{equation}
\begin{split}
\frac{-1 + c_0 + 2c_1\,X + c_2\,X^2}2 
&= \frac{(X-\xi)\bigl(\lambda(X-\eta)+\xi-\eta\bigr)}{(\xi-\eta)^2}\,,\\
\frac{-1 - c_0 - 2c_1\,Y - c_2\,Y^2}2
&= \frac{(\eta-Y)\bigl(\lambda(Y-\xi)+\xi-\eta\bigr)}{(\xi-\eta)^2}\,.
\end{split}
\end{equation}
Moreover, $X$ and~$Y$ are required to satisfy~$Y<\eta<\xi<X$.  
If~$\lambda=0$, this is the only restriction needed to make the right hand
sides positive, so the notation~$D_{\xi,\eta,0}$ will denote the
quarter-plane~$Y<\eta<\xi<X$.

If~$0<\lambda<1$, then~$Y$ is required to lie in
the interval
\begin{equation}
\frac\eta\lambda+\xi\left(1-\frac1\lambda\right)<Y<\eta\,
\end{equation}
so the notation~$D_{\xi,\eta,\lambda}$ will denote the corresponding 
open semi-infinite horizontal strip in the~$XY$-plane.

Finally, if~$-1<\lambda<0$, then~$X$ is required to lie in the interval
\begin{equation}
\xi<X<-\frac\xi\lambda+\eta\left(1-\frac1\lambda\right)\,,
\end{equation}
so the notation~$D_{\xi,\eta,\lambda}$ will denote the corresponding 
open semi-infinite vertical strip in the~$XY$-plane.

In the other direction, note that, if~$(X_0,Y_0)$ is any point 
in the~$XY$-plane with~$X_0>Y_0$, then the inequalities
\begin{equation}
-1 + c_0 + 2c_1\,X_0 + c_2\,{X_0}^2 > 0,\qquad
-1 - c_0 - 2c_1\,Y_0 - c_2\,{Y_0}^2 > 0
\end{equation}
define a non-empty open wedge in~$c_0c_1c_2$-space.  
Thus, for any point~$(X_0,Y_0)$ in the half-plane, 
there are choices of~$c_0$, $c_1$, and~$c_2$ so that the 
corresponding functions~$a$ and~$b$ are positive on a neighborhood
of~$(X_0,Y_0)$.  

\subsubsection{The connection forms}\label{sssec: connection forms I}
The connection structure forms have the expressions
\begin{equation}\label{eq: connection omegas in XY-terms I}
\begin{split}
\omega_{31} 
   &= \frac{\theta_1}{\sqrt{a}} 
 = \frac{(\xi-\eta)\,\d X}
      {(X-Y)\sqrt{(X-\xi)\bigl(\lambda(X-\eta)+\xi-\eta\bigr)\strut}} \,,\\
\omega_{32} 
   &=\frac{\theta_2}{\sqrt{b}}
 =\frac{(\xi-\eta)\,\d Y}
      {(Y-X)\sqrt{(\eta-Y)\bigl(\lambda(Y-\xi)+\xi-\eta\bigr)\strut}}\,,\\
\omega_{12} &= -\frac{\sqrt{b}}{\sqrt{a}}\,\theta_1
               +\frac{\sqrt{a}}{\sqrt{b}}\,\theta_2\\
    &= \frac{-1}{X{-}Y}\left(
\sqrt{\frac{(\eta-Y)\bigl(\lambda(Y-\xi)+\xi-\eta\bigr)}
           {(X-\xi)\bigl(\lambda(X-\eta)+\xi-\eta\bigr)}}\,\d X\right.\\
&\qquad\qquad\qquad\qquad
+ \left.\sqrt{\frac{(X-\xi)\bigl(\lambda(X-\eta)+\xi-\eta\bigr)}
                   {(\eta-Y)\bigl(\lambda(Y-\xi)+\xi-\eta\bigr)}}\,\d Y
\right).
\end{split}
\end{equation}

Note that the formulae for~$\omega_{31}$, $\omega_{32}$, and~$\omega_{12}$
do not explicitly involve the functions~$U$ and~$V$.%
\footnote{Of course, these
functions were used in the definition of the forms~$\theta_1$ and~$\theta_2$
and thus in the definition of the principal coordinatization~$(X,Y)$.}
In particular, the structure equations 
for these forms are satisfied on~$D_{\xi,\eta,\lambda}$.

\subsubsection{Euler linearization}\label{sssec: euler linear I}
As was seen in the general case, the structure equations 
for~$\omega_1$ and~$\omega_2$ (which are all that remains)
simplify to the linear system
\begin{equation}\label{eq: uv linear system I}
\frac{\partial U}{\partial Y} = -\frac{U-V}{X-Y}\,,\qquad\qquad 
\frac{\partial V}{\partial X} = -\frac{U-V}{X-Y}\,.
\end{equation}  
The general solution described in~\eqref{eq: UV explicit}
now simplifies to
\begin{equation}\label{eq: uv gen soln I}
\begin{split}
U(X,Y) &= f(X) - g(Y) - (X{-}Y)\,f'(X)\,,\\
V(X,Y) &= f(X) - g(Y) - (X{-}Y)\,g'(Y)\,,
\end{split}
\end{equation}
where~$f$ and~$g$ are arbitrary functions of a single variable, subject
only to the conditions that they be chosen on their respective $X$-domain 
and~$Y$-domain so that~$U$ and~$V$ are nonzero and nowhere equal functions 
on~$D$ (which is the product of the~$X$-domain and the $Y$-domain).
The functions~$f$ and~$g$ that give rise to~$U$ and~$V$ are not unique.  
In fact, for any constants~$m_0$ and~$m_1$ one can add~$m_0+m_1\,X$ to~$f(X)$ 
and~$m_0+m_1\,Y$ to~$g(Y)$ without changing~$U$ and~$V$.  However, this
is the only indeterminacy in the formulae.

The result is the following general formulae for~$\omega_1$ and~$\omega_2$:
\begin{equation}\label{eq: trans omegas in XY-terms I}
\begin{split}
\omega_1 &= \frac{\theta_1}{A\sqrt{a}}
      = \frac{(\xi-\eta)\,\bigl(f(X) - g(Y) - (X{-}Y)\,f'(X)\bigr)\,\d X}
      {(X-Y)\sqrt{(X-\xi)\bigl(\lambda(X-\eta)+\xi-\eta\bigr)\strut}} \,,\\
\omega_2 &= \frac{\theta_2}{B\sqrt{b}}
      = \frac{(\xi-\eta)\,\bigl(f(X) - g(Y) - (X{-}Y)\,g'(Y)\bigr)\,\d Y}
      {(Y-X)\sqrt{(\eta-Y)\bigl(\lambda(Y-\xi)+\xi-\eta\bigr)\strut}}\,,\,,\\
\end{split}
\end{equation}

Note that the forms~$\omega_1,\omega_2,\omega_{31},\omega_{32},\omega_{12}$
as defined in~\eqref{eq: connection omegas in XY-terms I} 
and~\eqref{eq: trans omegas in XY-terms I} satisfy the structure equations
even at points where~$U$ or~$V$ vanish or where~$U=V$ (i.e., 
where~$f'(X) = g'(Y)$).  Consequently,
Bonnet's theorem applies and there exist mappings~$\xb:D\to\E{3}$
and~$(\eb_1,\eb_2,\eb_3):D\to\Or(3)$ whose associated structure
forms are the given ones.  As long as~$U$ and~$V$ are nonzero, the
map~$\xb$ will be an immersion, it is just that
this immersion will have umbilic points where~$f'(X)=g'(Y)$.

\subsubsection{A Weierstra{\ss}-type formula}\label{sssec: weierstrass I}
Now, it is not necessary to rely on Bonnet's theorem to generate the
mapping~$\xb$.  In fact, this can be reduced to quadratures, as will now
be demonstrated.

In the first place, 
finding a frame field~$\eb = (\eb_1,\eb_2,\eb_3):D\to\Or(3)$ so
that~$\d \eb_i = \eb_j\,\omega_{ji}$ (where~$\omega_{ji} = -\omega_{ij}$)
can be done as follows.  One starts with the formula
\begin{equation}\label{eq: classical e3 I}
\eb_3 = \left[
\,\begin{matrix}
\ds\frac{2\sqrt{(X-\xi)\bigl(\lambda(X-\eta)+\xi-\eta\bigr)\strut}}
     {(1+\lambda)(X-Y)}\\[15pt]
\ds\frac{2\sqrt{(\eta-Y)\bigl(\lambda(Y-\xi)+\xi-\eta\bigr)\strut}}
     {(1-\lambda)(X-Y)}\\[15pt]
\ds\frac{2(\xi-\eta)-2\lambda(\xi+\eta)-(1-\lambda)^2X+(1+\lambda)^2Y}
     {(1-\lambda^2)(X-Y)}
\end{matrix}\,\right].
\end{equation}

Computation yields
\begin{equation}
\begin{split}
\frac{\partial\eb_3}{\partial X}\cdot\frac{\partial\eb_3}{\partial X}
&=-\frac{(\xi-\eta)^2}
        {(X-Y)^2(X-\xi)\bigl(\lambda(X-\eta)+\xi-\eta\bigr)}\,,\\
\frac{\partial\eb_3}{\partial X}\cdot\frac{\partial\eb_3}{\partial Y}
&=0,\\
\frac{\partial\eb_3}{\partial Y}\cdot\frac{\partial\eb_3}{\partial Y} 
&= \frac{(\xi-\eta)^2}
        {(X-Y)^2(\eta-Y)\bigl(\lambda(Y-\xi)+\xi-\eta\bigr)}\,,
\end{split}
\end{equation}
so setting
\begin{equation}
\eb_1 = -\left|\frac{\partial\eb_3}{\partial X}\right|^{-1}\,
         \frac{\partial\eb_3}{\partial X}\,,\qquad\qquad
\eb_2 = -\left|\frac{\partial\eb_3}{\partial Y}\right|^{-1}\,
         \frac{\partial\eb_3}{\partial Y}\,,
\end{equation}
yields the desired structure equation
\begin{equation}
\d\eb_3 = \frac{\partial\eb_3}{\partial X}\,\d X
        + \frac{\partial\eb_3}{\partial Y}\,\d Y
  = {}-\eb_1\,\omega_{31}-\eb_2\,\omega_{32}\,,
\end{equation}
where~$\omega_{31}$ and~$\omega_{32}$ are as defined 
in~\eqref{eq: connection omegas in XY-terms I}.

The remaining structure equations for~$\d\eb_1$ and~$\d\eb_2$
are easily verified, so that this does, in fact, integrate the
structure equations for the frame field~$(\eb_1,\eb_2,\eb_3)$.

\begin{remark}[The Gau{\ss} image of the principal net]
\label{rem: Gauss prin curves I}
Note that, since~$\mu=0$, both families of principal
curves are mapped to arcs of circles under~$\eb_3$.
These two families of image circles are evidently orthogonal.

In particular, for any realization~$\xb$ of
a Type~I shape operator~$S$ with~$\mu=0$, 
the Gau{\ss} image of the net of principal curves of~$\xb$ 
is a net of orthogonal circle foliations on the $2$-sphere.
\end{remark}

Now, note that~$\eb_3$ satisfies the (vector-valued) Euler equation 
\begin{equation}\label{eq: e3 Euler I}
\frac{\partial^2 \eb_3}{\partial X\partial Y} 
- \frac{1}{X-Y}\,\frac{\partial \eb_3}{\partial X}
+ \frac{1}{X-Y}\,\frac{\partial \eb_3}{\partial Y}
 = 0,
\end{equation}
which is easily established by direct computation.

Finally, the immersion~$\xb:D\to\E{3}$ satisfies the structure equation
\begin{equation}\label{eq: dx equation I}
\begin{split}
\d\xb &= \eb_1\,\omega_1 + \eb_2\,\omega_2 
      = U(X,Y)\,\eb_1\,\omega_{31} + V(X,Y)\,\eb_2\,\omega_{32}\\
      &= {}-U(X,Y)\,\frac{\partial\eb_3}{\partial X}\,\d X
           -V(X,Y)\,\frac{\partial\eb_3}{\partial Y}\,\d Y.
\end{split}
\end{equation}
The fact that~$\eb_3$ satisfies~\eqref{eq: e3 Euler I} 
while~$U$ and~$V$ satisfy~\eqref{eq: uv linear system I} implies
the identity
\begin{equation}\label{eq: e3 compat I}
\frac{\partial\hfil}{\partial Y}
      \left(U(X,Y)\,\frac{\partial\eb_3}{\partial X}\right)
= 
\frac{\partial\hfil}{\partial X}
      \left(V(X,Y)\,\frac{\partial\eb_3}{\partial Y}\right),
\end{equation}
implying that the vector-valued differential form on the right hand
side of~\eqref{eq: dx equation I} is indeed a closed $1$-form, so
that~$\xb$ can be recovered by quadrature:
\begin{equation}\label{eq: x quadrature I}
\xb = - \int U(X,Y)\,\frac{\partial\eb_3}{\partial X}\,\d X
           +V(X,Y)\,\frac{\partial\eb_3}{\partial Y}\,\d Y\,.
\end{equation}

\subsubsection{Unfoldings}\label{sssec: unfoldings I}
If a solution~$(U,V)$ to~\eqref{eq: uv linear system I}
is defined on a neighborhood of the closure of a 
domain~$D_{\xi,\eta,\lambda}$, then the mapping~$\xb$ can be
continued beyond the edges of~$D_{\xi,\eta,\lambda}$ in a
natural way.  

For example, when~$\lambda=0$, consider the mapping from~$\bbR^2$
to the closure of~$D_{\xi,\eta,0}$ defined by the formulae
\begin{equation}\label{eq: lambda=0 unfolding}
\begin{split}
X &= \xi  + x^2\,,\\
Y &= \eta - y^2\,.
\end{split}
\end{equation}
Using this mapping to pull back~$\eb_3$ (assuming that~$\lambda=0$,
of course), the formula for~$\eb_3$ is resolvable to
\begin{equation}\label{eq: lam=0 unfolded e3 I}
\eb_3(x,y) = \left[
\,\begin{matrix}
\ds\frac{2x\,\sqrt{\xi{-}\eta}}{\bigl((\xi{-}\eta)+x^2+y^2\bigr)}\\[10pt]
\ds\frac{2y\,\sqrt{\xi{-}\eta}}{\bigl((\xi{-}\eta)+x^2+y^2\bigr)}\\[10pt]
\ds \frac{(\xi{-}\eta)-x^2-y^2}{\bigl((\xi{-}\eta)+x^2+y^2\bigr)}
\end{matrix}\,\right],
\end{equation}
and the reader will notice that this is a conformal embedding of the 
$xy$-plane onto the punctured sphere.  The differential formula
for~$\xb$ as a function of~$x$ and~$y$ then becomes
\begin{equation}\label{eq: x unfolded lambda=0}
\d\xb = U(\xi{+}x^2,\eta{-}y^2)\,\frac{\partial\eb_3}{\partial x}\,\d x
       +V(\xi{+}x^2,\eta{-}y^2)\,\frac{\partial\eb_3}{\partial y}\,\d y\,,
\end{equation}
and the $1$-form on the right hand side of this equation will be smooth
and closed as long as~$U$ and~$V$ satisfy~\eqref{eq: uv linear system I}
and are smooth on a domain containing the closure of~$D_{\xi,\eta,0}$.
Moreover, it will be an immersion on the set where~$U$ and~$V$ are 
nonzero.

There are similar unfoldings when~$\lambda\not=0$.  When~$\lambda>0$,
one uses the formulae
\begin{equation}\label{eq: lambda>0 unfolding}
\begin{split}
X &= \xi  + x^2\,,\\
Y &= \eta\,\cos^2y 
+\left(\frac\eta\lambda +\xi\left(1-\frac1\lambda\right)\right)\,\sin^2y\,,
\end{split}
\end{equation}
while, when~$\lambda<0$, one uses the formulae
\begin{equation}\label{eq: lambda<0 unfolding}
\begin{split}
X &= \xi\,\cos^2y 
+\left(-\frac\xi\lambda +\eta\left(1-\frac1\lambda\right)\right)\,\sin^2y\,,\\
Y &= \eta - y^2\,,
\end{split}
\end{equation}
to obtain mappings defined on a smooth cylinder.  In either case,~$\eb_3$
becomes a conformal embedding of the cylinder onto the twice-punctured unit
$2$-sphere and the formula analogous to~\eqref{eq: x unfolded lambda=0}
defines (up to a quadrature) a smooth mapping~$\xb$ from the cylinder
to~$\E3$.  

\subsubsection{Examples}\label{sssec: examples I}

Some examples will now be considered.

\begin{example}[Linear solutions]\label{ex: linear type I}
Consider the global linear solution to~\eqref{eq: uv linear system I}
\begin{equation}\label{eq: linear solns I}
U = Y\,,\qquad\qquad\qquad V = X\,,
\end{equation} 
whose domain is the entire half-plane~$Y<X$.  Using the formulae
above, one finds the corresponding mapping on~$D_{\xi,\eta,\lambda}$
to be given by the formulae
\begin{equation}\label{eq: x linear solns I}
\xb(X,Y) = \left[
\,\begin{matrix}
\ds\frac{-2Y\,\sqrt{(X-\xi)\bigl(\lambda(X-\eta)+\xi-\eta\bigr)\strut}}
     {(1+\lambda)(X-Y)}\\[15pt]
\ds\frac{-2X\,\sqrt{(\eta-Y)\bigl(\lambda(Y-\xi)+\xi-\eta\bigr)\strut}}
     {(1-\lambda)(X-Y)}\\[15pt]
\ds\frac{(X+Y)\bigl(\lambda(\xi+\eta)-\xi+\eta\bigr)
       -4\lambda\,XY}
     {(1-\lambda^2)(X-Y)}
\end{matrix}\,\right].
\end{equation}
This is an immersion away from the axis rays~$X=0$ and~$Y=0$.  
The shape operator of all of these immersions is
\begin{equation}\label{eq: S linear solns I}
S = \frac1Y\, \frac{\partial\hfil}{\partial X}\otimes \d X
   +\frac1X\, \frac{\partial\hfil}{\partial Y}\otimes \d Y\,.
\end{equation}

Note that, since the domain of this solution contains the closures 
of all of the~$D_{\xi,\eta,\lambda}$, it follows that the various
unfoldings allow one to continue~$\xb$ past the edges in all
cases as a smooth mapping.  This mapping will be an immersion as
long as the closure of~$D_{\xi,\eta,\lambda}$ does not meet either
axis.

The image surface is visibly algebraic.  However, it does not 
appear to be easy to recognize as a classical surface.  Indeed,
for generic values of~$\xi$, $\eta$, and~$\lambda$, it appears to
be of degree~$8$.

\end{example}

\begin{example}[An umbilic of index zero]\label{ex: index zero type I}
Consider the polynomial solution~$(U,V)$ generated
by the formulae~\eqref{eq: uv gen soln I} when one takes
\begin{equation}
f(X) = 1+(X-1)^3\,,\qquad g(Y) = 1-(Y+1)^3.
\end{equation}
Since~$U=V$ only where~$f'(X) = g'(Y)$, one sees that this happens
only at the point~$(X,Y) = (1,-1)$.  Moreover, $U$ and~$V$ are
positive at this point.  
It follows that there is a neighborhood of~$(1,-1)$ in the~$XY$-plane
on which the shape operator
\begin{equation}\label{eq: S umbilic index zero solns I}
S = \frac1{U(X,Y)}\,\frac{\partial\hfil}{\partial X}\otimes \d X
   +\frac1{V(X,Y)}\,\frac{\partial\hfil}{\partial Y}\otimes \d Y\,
\end{equation}
has a $3$-parameter family of non-congruent realizations.  All of
these realizations have an isolated umbilic point and the construction
shows that this umbilic point is of index zero. (After all, 
the net of principal curves is non-singular near this point.)

Whether such an example exists globally is an interesting question.
\end{example}

\begin{example}[A minimal surface]\label{ex: minimal surface I}
It is not difficult to see that, up to a constant multiple, 
the only solution~$(U,V)$ to~\eqref{eq: uv gen soln I} 
that satisfies~$U+V=0$ is 
\begin{equation}
U(X,Y) = (X-Y)^2\,,\qquad V(X,Y) = -(X-Y)^2.
\end{equation}
(Note that this solution is invariant under the simultaneous
translation of~$X$ and~$Y$ and simply scales under the simultaneous
dilation.)  Note that~$U+V = (A+B)/(AB)$, so that~$U+V=0$ if
and only if~$A+B=0$, i.e., if the realization is a minimal surface.%
\footnote{Adding a constant~$c$ to each of~$U$ and~$V$ yields
another solution that satisfies~$U+V=2c$.  Since this is
equivalent to~$(A+B)/(AB)= 2c$, i.e., to~$H = cK$ (where~$H$
and~$K$ are the mean and Gau{\ss} curvatures respectively),
this gives a more general class of Weingarten surfaces that
admit a 3-parameter family of deformations.  It is not difficult
to show that there is, in fact, a two parameter family of 
Weingarten relations that have deformable examples of this type.}

Thus, this solution gives, up to constant multiples, the only
shape operator of Type~I whose realizations are minimal surfaces.%
\footnote{The argument given only applies to the~$\mu=0$ case, but
the reader will have no difficulty checking that when~$\mu\not=0$,
there are no solutions with~$U+V=0$.}

When~$\lambda$ is nonzero, the integral that gives~$\xb$ for this
solution is rather complicated, so it will not be written out here.
Instead, I will simply note that the integrals in the case~$\lambda=0$
give rise to Enneper's surface (up to translation and scale), 
as the reader can easily verify. 
\end{example}

\subsection{Operators of Type~II}\label{ssec: Type II}
I now want to consider the operators of Type~II and explain how the
equations that define them can be integrated to a linear system. 

\subsubsection{Natural principal coordinates}
\label{sssec: nat prin coords II}
Recall the definition of the $1$-forms
\begin{equation}
\theta_1  = \frac{A\,B_x}{B(B-A)}\,\d x\,,\qquad\qquad
\theta_2 =  \frac{B\,A_y}{A(A-B)}\,\d y\,.
\end{equation}
and that the equations~\eqref{eq: second system} can be expressed as 
\begin{equation}
\d \theta_1 =  -2\,\theta_1 \w \theta_2\,,\qquad\qquad
\d \theta_2 =  -2\,\theta_2 \w \theta_1\,.
\end{equation}
As a consequence, there exist functions~$X$ and~$Y$ on~$D$ with~$X>Y$
and for which 
\begin{equation}\label{eq: theta II in XY-terms}
\theta_1 = \frac{\d X}{2(Y-X)}\,,\qquad\qquad
\theta_2 = \frac{\d Y}{2(X-Y)}\,.
\end{equation}

In fact, these two functions can be found by quadrature as follows:
First, the structure equations imply that the $1$-form~$\theta_1+\theta_2$
is closed, so, by quadrature, one can find a function~$Z>0$ on~$D$
so that~$\d\bigl(\log Z\bigr) = \d Z/Z = -2(\theta_1+\theta_2)$.
The structure equations now imply that the~$1$-form $-2Z\,\theta_1$
is closed, so, again, by quadrature, one can find a function~$X$
so that~$\d X = -2Z\,\theta_1$.  Setting~$Y = X-Z$ then gives the
desired remaining function.

Note that~$X$ and~$Y$ are unique up to a replacement 
of the form~$(X,Y)\mapsto (\lambda\,X+\tau,\lambda\,Y+\tau)$
where~$\lambda>0$ and~$\tau$ are constants.  
The map~$(X,Y):D\to\bbR^2$ is a principal coordinatization of~$D$ 
and embeds~$D$ as a rectangle in the open half-plane~$X>Y$.
As before, these will be referred to as \emph{natural principal
coordinates}.

\subsubsection{Integrals of the Frobenius system}\label{sssec: Frob ints II}
Now, the Frobenius system for the functions~$a$,~$b$, and~$p$
simplifies to
\begin{equation}\label{eq: abp Frobenius II}
\begin{split}
\d a & =    \phantom{4(2b+p+1)}\hbox to 0pt{\hss$(2a-p+1)$}\,\theta_1
         +  \phantom{4(2a-p+1)}\hbox to 0pt{\hss$6a\,$}\,\theta_2\,,\\
\d b & =    \phantom{4(2b+p+1)}\hbox to 0pt{\hss$6b\,$}\,\theta_1 
         +  \phantom{4(2a-p+1)}\hbox to 0pt{\hss$(2b+p+1)$}\,\theta_2\,,\\
\d p & = 4(2b+p+1)\,\theta_1 - 4(2a-p+1)\,\theta_2\,.\\
\end{split}
\end{equation}

By standard integration techniques, 
there exist constants~$c_0$, $c_1$, and $c_2$ so that
\begin{equation}\label{eq: abp in XY-terms II}
\begin{split}
a &= \frac{X^3+3c_2\,X^2+3c_1\,X+c_0}{(Y-X)^3}\,,\\
b &= \frac{Y^3+3c_2\,Y^2+3c_1\,Y+c_0}{(X-Y)^3}\,,\\
p &= \frac{Y^3-3Y^2X-3YX^2+X^3-12c_2\,XY-6c_1\,(X+Y)-4c_0}{(Y-X)^3}\,.
\end{split}
\end{equation}
Conversely, for any constants~$c_0$, $c_1$, and $c_2$, 
the formulae~\eqref{eq: abp in XY-terms II} give expressions 
for~$a$, $b$, and~$p$ that satisfy~\eqref{eq: abp Frobenius II}.
The reader should bear in mind, however, that the conditions~$a>0$
and~$b>0$ must still be imposed.  These inequalities, together with
the requirement~$X>Y$, impose inequalities on~$c_0$,~$c_1$,~and~$c_2$
that amount to requiring that the 
polynomial~$c(t) = t^3 + 3c_2\,t^2+3c_1\,t+c_0$ have three distinct
real roots, say,
\begin{equation}\label{eq: c factored II}
c(t) = (t-\lambda_1)(t-\lambda_2)(t-\lambda_3), \qquad\qquad
\lambda_1 < \lambda_2 < \lambda_3
\end{equation}
If these requirements are met, then~$a$ and~$b$ are positive 
on the open rectangle~$D_\lambda$ defined by the inequalities
\begin{equation}\label{eq: XY limits II}
\lambda_1 < Y < \lambda_2 < X < \lambda_3\,.
\end{equation}
Note that one corner of~$D_\lambda$ lies on the line~$X=Y$, 
namely the point~$(\lambda_2,\lambda_2)$.

In the other direction, note that, if~$(X_0,Y_0)$ is any point 
in the~$XY$-plane with~$X_0>Y_0$, then the inequalities
\begin{equation}
{X_0}^3+3c_2\,{X_0}^2+3c_1\,X_0+c_0 < 0,\qquad
{Y_0}^3+3c_2\,{Y_0}^2+3c_1\,Y_0+c_0 > 0
\end{equation}
define a non-empty open wedge in~$c_0c_1c_2$-space.  
Thus, for any point~$(X_0,Y_0)$ in the half-plane, 
there are choices of~$c_0$, $c_1$, and~$c_2$ so that the 
corresponding functions~$a$ and~$b$ are positive on a neighborhood
of~$(X_0,Y_0)$. 

\subsubsection{The connection forms}\label{sssec: connection forms II}
Now, in terms of~$\theta_1$ and~$\theta_2$, 
the connection structure forms have the expressions
\begin{equation}\label{eq: connection omegas in XY-terms II}
\begin{split}
\omega_{31} 
   &= \frac{\theta_1}{\sqrt{a}} 
     = \frac{\sqrt{(X-Y)}\,\d X}
         {2\sqrt{(X-\lambda_1)(X-\lambda_2)(\lambda_3-X)\strut}} \,,\\
\omega_{32} 
   &=\frac{\theta_2}{\sqrt{b}}
    =\frac{\sqrt{(X-Y)}\,\d Y}
         {2\sqrt{(Y-\lambda_1)(\lambda_2-Y)(\lambda_3-Y)\strut}}\,,\\
\omega_{12} &= -\frac{\sqrt{b}}{\sqrt{a}}\,\theta_1
               +\frac{\sqrt{a}}{\sqrt{b}}\,\theta_2\\
    &= \frac{1}{2(X{-}Y)}\left(
  \sqrt{\frac{(Y-\lambda_1)(\lambda_2-Y)(\lambda_3-Y)}
             {(X-\lambda_1)(X-\lambda_2)(\lambda_3-X)}}\,\d X\right.\\
&\qquad\qquad\qquad\qquad
+ \left.\sqrt{\frac{(X-\lambda_1)(X-\lambda_2)(\lambda_3-X)}
                   {(Y-\lambda_1)(\lambda_2-Y)(\lambda_3-Y)}}\,\d Y
\right).
\end{split}
\end{equation}

Note that the formulae for~$\omega_{31}$, $\omega_{32}$, and~$\omega_{12}$
do not explicitly involve the functions~$A$ and~$B$.%
\footnote{Of course, these
functions were used in the definition of the forms~$\theta_1$ and~$\theta_2$
and thus in the definition of the principal coordinatization~$(X,Y)$.}
In particular, the structure equations for these forms are satisfied 
on~$D_\lambda$.

\subsubsection{Euler linearization}\label{sssec: euler linear II}
Writing~$U=1/A$ and~$V=1/B$, the structure equations 
for~$\omega_1$ and~$\omega_2$ (which are all that remains)
simplify to the linear system
\begin{equation}\label{eq: uv linear system II}
\frac{\partial U}{\partial Y} = \frac{U-V}{2(X-Y)}\,,\qquad\qquad 
\frac{\partial V}{\partial X} = \frac{U-V}{2(X-Y)}\,.
\end{equation} 
However, this linear system is \emph{not} integrable by 
the method of Darboux,  so its general solution cannot
be expressed in a closed form similar to~\eqref{eq: uv gen soln I}.
One can express the system as a single hyperbolic equation by 
introducing a potential~$\Phi(X,Y)$ so that~$U=\Phi_X$ and~$V=\Phi_Y$. 
Then~$\Phi$ satisfies the so-called Euler equation
\begin{equation}\label{eq: Phi Euler}
\frac{\partial^2 \Phi}{\partial X\partial Y} 
- \frac{\textstyle{\frac12}}{X-Y}\,\frac{\partial \Phi}{\partial X}
+ \frac{\textstyle{\frac12}}{X-Y}\,\frac{\partial \Phi}{\partial Y}
 = 0.
\end{equation}
While there is no closed-form solution to this equation, Poisson
has given the following integral formula for the general solution
\begin{equation}\label{eq: Poisson formula}
\begin{split}
\Phi(X,Y) 
&= \int_Y^X\frac{\phi(\xi)}{\sqrt{(X-\xi)(\xi-Y)}}\,d\xi\\
&\qquad\qquad + \int_Y^X
  \frac{\psi(\xi)}{\sqrt{(X-\xi)(\xi-Y)}}
   \,\log\left(\frac{(X-Y)}{(X-\xi)(\xi-Y)}\right)\,d\xi\,,
\end{split}
\end{equation}
where~$\phi$ and~$\psi$ are arbitrary functions of a single variable.
If~$\phi$ and~$\psi$ are defined on an interval~$(a,b)$, then the
solution~$\Phi$ is defined on the triangle~$a<Y<X<b$.  (Of course,
$a=+\infty$ and/or~$b=-\infty$ are allowable values.)

These formulae give the general solution (in natural principal coordinates)
to the system~\eqref{eq: second system}.  Thus, every operator~$S$ of Type~II
can be generated by this procedure.

The result is the following general formulae for~$\omega_1$ and~$\omega_2$:
\begin{equation}\label{eq: trans omegas in XY-terms II}
\begin{split}
\omega_1 &= \frac{\theta_1}{A\sqrt{a}}
      = \frac{\Phi_X(X,Y)\,\sqrt{(X-Y)}\,\d X}
              {2\sqrt{(X-\lambda_1)(X-\lambda_2)(\lambda_3-X)}} \,,\\
\omega_2 &= \frac{\theta_2}{B\sqrt{b}}
      = \frac{\Phi_Y(X,Y)\,\sqrt{(X-Y)}\,\d Y}
              {2\sqrt{(Y-\lambda_1)(\lambda_2-Y)(\lambda_3-Y)}}\,,\,,\\
\end{split}
\end{equation}

Note that the forms~$\omega_1,\omega_2,\omega_{31},\omega_{32},\omega_{12}$
as defined in~\eqref{eq: connection omegas in XY-terms II} 
and~\eqref{eq: trans omegas in XY-terms II} satisfy the structure equations
even at points where~$\Phi_X$, or~$\Phi_Y$ vanish
or where~$\Phi_X=\Phi_Y$, as long as~$\Phi$
satisfies~\eqref{eq: Phi Euler}.   Consequently, Bonnet's theorem
applies and there exist mappings~$\xb: D\to \E3$ 
and~$(\eb_1,\eb_2,\eb_3): D\to\Or(3)$ whose associated structure forms
are the given ones.  As long as $\Phi_X$ and~$\Phi_Y$ are nonzero,
the map~$\xb$ will be an immersion, it is just that this immersion
will have umbilic points where~$\Phi_X=\Phi_Y$.

\subsubsection{A Weierstra{\ss}-type formula}\label{sssec: weierstrass II}
Now, it is not necessary to rely on Bonnet's theorem to generate the
mapping~$\xb$.  In fact, this can be reduced to quadratures, as will
now be demonstrated.

In the first place, finding 
a frame field~$\eb = (\eb_1,\eb_2,\eb_3):D\to \SO(3)$ so
that~$\d \eb_i = \eb_j\,\omega_{ji}$ (where~$\omega_{ji} = -\omega_{ij}$)
is easily done.  One starts with the classical formula
\begin{equation}\label{eq: classical e3}
\eb_3 = \left[\,\,\begin{matrix}
\ds\sqrt{\frac{(X{-}\lambda_1)(Y{-}\lambda_1)}
           {(\lambda_2{-}\lambda_1)(\lambda_3{-}\lambda_1)}}\\[15pt]
\ds\sqrt{\frac{(X{-}\lambda_2)(\lambda_2{-}Y)}
           {(\lambda_2{-}\lambda_1)(\lambda_3{-}\lambda_2)}}\\[15pt]
\ds\sqrt{\frac{(\lambda_3{-}X)(\lambda_3{-}Y)}
           {(\lambda_3{-}\lambda_1)(\lambda_3{-}\lambda_2)}}
\end{matrix}\,\,\right]\,.
\end{equation}

Note that~$\eb_3\cdot\eb_3=1$.  In fact,
$\eb_3$ maps~$D_\lambda$ diffeomorphically onto the positive
orthant of the $2$-sphere and extends continuously to the boundary
of~$D_\lambda$ as a homeomorphism from the closure of~$D_\lambda$
to the closure of the positive orthant.   

One easily computes that
\begin{equation}
\begin{split}
\frac{\partial\eb_3}{\partial X}\cdot\frac{\partial\eb_3}{\partial X}
&= -\frac{X-Y}{4(X-\lambda_1)(X-\lambda_2)(X-\lambda_3)}\,,\\
\frac{\partial\eb_3}{\partial X}\cdot\frac{\partial\eb_3}{\partial Y}
&=0,\\
\frac{\partial\eb_3}{\partial Y}\cdot\frac{\partial\eb_3}{\partial Y} 
&= \frac{X-Y}{4(Y-\lambda_1)(Y-\lambda_2)(Y-\lambda_3)}\,,
\end{split}
\end{equation}
so setting
\begin{equation}
\eb_1 = -\left|\frac{\partial\eb_3}{\partial X}\right|^{-1}\,
         \frac{\partial\eb_3}{\partial X}\,,\qquad\qquad
\eb_2 = -\left|\frac{\partial\eb_3}{\partial Y}\right|^{-1}\,
         \frac{\partial\eb_3}{\partial Y}\,,
\end{equation}
yields the desired structure equation
\begin{equation}
\d\eb_3 = \frac{\partial\eb_3}{\partial X}\,\d X
        + \frac{\partial\eb_3}{\partial Y}\,\d Y
  = {}-\eb_1\,\omega_{31}-\eb_2\,\omega_{32}\,,
\end{equation}
where~$\omega_{31}$ and~$\omega_{32}$ are as defined 
in~\eqref{eq: connection omegas in XY-terms II}.

The remaining structure equations for~$\d\eb_1$ and~$\d\eb_2$ are
easily verified, so that this does, in fact, integrate the equations
for the frame field~$(\eb_1,\eb_2,\eb_3)$.

One very interesting consequence of having this explicit 
form of~$\eb_3$ is the following result due to Finikoff and Gambier:

\begin{proposition}\label{prop:  Gauss prin curves II}
For any realization~$\xb$ of a Type~II shape operator~$S$, 
the Gau{\ss} image of the net of principal curves of~$\xb$ 
is an orthogonal net of confocal spherical ellipses.
\end{proposition}

\begin{proof} 
Write
\begin{equation}
\eb_3(X,Y) 
= \left[\begin{matrix} u(X,Y)\\ v(X,Y)\\ w(X,Y)\end{matrix}\right]
\end{equation} 
and note that these components satisfy the equations
\begin{equation}\label{eq: confocal relations}
\begin{split}
\frac{u^2}{X-\lambda_1}+\frac{v^2}{X-\lambda_2}+\frac{w^2}{X-\lambda_3}&=0,\\
\frac{u^2}{Y-\lambda_1}+\frac{v^2}{Y-\lambda_2}+\frac{w^2}{Y-\lambda_3}&=0.
\end{split}
\end{equation}
This is the content of the proposition.
\end{proof}

Now, note that~$\eb_3$ satisfies the (vector-valued) Euler equation 
that is dual to the Euler equation~\eqref{eq: Phi Euler} satisfied by~$\Phi$
\begin{equation}\label{eq: e3 Euler II}
\frac{\partial^2 \eb_3}{\partial X\partial Y} 
- \frac{\textstyle{-\frac12}}{X-Y}\,\frac{\partial \eb_3}{\partial X}
+ \frac{\textstyle{-\frac12}}{X-Y}\,\frac{\partial \eb_3}{\partial Y}
 = 0,
\end{equation}
which is easily established by direct computation.

Finally, the immersion~$\xb:D\to\E{3}$ satisfies the structure equation
\begin{equation}\label{eq: dx equation II}
\begin{split}
\d\xb &= \eb_1\,\omega_1 + \eb_2\,\omega_2 
      = \Phi_X\,\eb_1\,\omega_{31} + \Phi_Y\,\eb_2\,\omega_{32}\\
      &= {}-\frac{\partial\Phi}{\partial X}\,
               \frac{\partial\eb_3}{\partial X}\,\d X
           -\frac{\partial\Phi}{\partial Y}\,
               \frac{\partial\eb_3}{\partial Y}\,\d Y.
\end{split}
\end{equation}
The fact that~$\eb_3$ and~$\Phi$ satisfy dual Euler equations implies
the identity
\begin{equation}
\frac{\partial\hfil}{\partial Y}
      \left(\frac{\partial\Phi}{\partial X}\,
               \frac{\partial\eb_3}{\partial X}\right)
= 
\frac{\partial\hfil}{\partial X}
      \left(\frac{\partial\Phi}{\partial Y}\,
               \frac{\partial\eb_3}{\partial Y}\right),
\end{equation}
implying that the vector-valued differential form on the right hand
side of~\eqref{eq: dx equation II} is indeed a closed $1$-form, so
that~$\xb$ can be recovered by quadrature:
\begin{equation}\label{eq: x quadrature II}
\xb = - \int \frac{\partial\Phi}{\partial X}\,
               \frac{\partial\eb_3}{\partial X}\,\d X
           +\frac{\partial\Phi}{\partial Y}\,
               \frac{\partial\eb_3}{\partial Y}\,\d Y\,.
\end{equation}

\subsubsection{Toral unfolding, spherical quotient}
\label{sssec: torus sphere II}
So far, no attention has been paid to the behavior of the immersion
near the edges of the rectangular domain~$D_\lambda$
defined by the inequalities~$\lambda_1<Y<\lambda_2<X<\lambda_3$.
It turns out, however, that one can often extend the immersion
to cover these edges.

Consider the mapping of the torus
\begin{equation}
T = \R{}/\left(2\pi\bbZ\right) \times \R{}/\left(2\pi\bbZ\right)
\end{equation}
into the closure of~$D_\lambda$ that is defined by the equations
\begin{equation}\label{eq: XY subs}
\begin{split}
Y = Y_\lambda(x,y) &= \lambda_1\,\sin^2y + \lambda_2\,\cos^2y\,,\\
X = X_\lambda(x,y) &= \lambda_2\,\cos^2x + \lambda_3\,\sin^2x \,.
\end{split}
\end{equation}
This mapping is not an immersion along the half-lattice lines, but
when one pulls~$\eb_3$ back to~$T$ by this mapping, the square roots
in the formula for~$\eb_3$ can be resolved into the form
\begin{equation}\label{eq: toral e3}
\eb_3(x,y) = \left[\,\,\begin{matrix}
\ds\cos y\,
\sqrt{\ds
\frac{(\lambda_3{-}\lambda_1)\,\sin^2x+(\lambda_2{-}\lambda_1)\,\cos^2x}
              {(\lambda_3{-}\lambda_1)}}\\[15pt]
\ds\sin x\,\,\sin y\\[15pt]
\ds\cos x\,
\sqrt{\ds
\frac{(\lambda_3{-}\lambda_1)\,\sin^2y+(\lambda_3{-}\lambda_2)\,\cos^2y}
              {(\lambda_3{-}\lambda_1)}}
\end{matrix}\,\,\right]\,,
\end{equation}
in which the expressions under the radicals are strictly positive 
on the torus.  Consequently, the map~$\eb_3:T\to S^2$ is smooth.  
Moreover, it is a submersion except at the four half-lattice points 
defined by the equations~$\sin x = \sin y = 0$.  
In fact, the map~$\eb_3$ is simply the quotient of~$T$ 
by the involution~$\tau:T\to T$ defined by $\tau(x,y)=(-x,-y)$, 
whose fixed points are the half-lattice points.  
Note also the symmetries
\begin{equation}\label{eq: e3 symmetries}
\eb_3(x{+}\pi, y) = R_u\,\eb_3(x, y) \,,\qquad\qquad
\eb_3(x, y{+}\pi) = R_w\,\eb_3(x, y)\,.
\end{equation}
where~$R_u$ and~$R_w$ are rotations by~$\pi$ about the $u$- and~$w$-axes,
respectively.  

If~$\Phi$ is a solution to~\eqref{eq: Phi Euler} that is smooth 
on a domain that contains the closure of~$D_\lambda$, 
then the formula~\eqref{eq: x quadrature II} 
can be lifted back to the torus in the form
\begin{equation}\label{eq: x torus quadrature II}
\begin{split}
\xb(x,y) &= - \int_{(0,0)}^{(x,y)} 
 \frac{\partial\Phi}{\partial X}
      \bigl(X_\lambda(\xi,\eta),Y_\lambda(\xi,\eta)\bigr)\,
       \frac{\partial\eb_3}{\partial x}(\xi,\eta)\,\d\xi\\
&\qquad\qquad\qquad\qquad
+\frac{\partial\Phi}{\partial Y}
      \bigl(X_\lambda(\xi,\eta),Y_\lambda(\xi,\eta)\bigr)\,
               \frac{\partial\eb_3}{\partial y}(\xi,\eta)\,\d\eta\,.
\end{split}
\end{equation}

Using the $\tau$-invariance of~$\eb_3$ 
and the symmetries~\eqref{eq: e3 symmetries},
is easy to show that the periods of the $1$-form integrand 
in~\eqref{eq: x torus quadrature II} must vanish.  Thus, the line
integral~\eqref{eq: x torus quadrature II} 
defines a smooth mapping~$\xb:T\to\E{3}$.

Moreover, since the integrand is invariant under~$\tau$,
it follows that~$\xb(x,y) = \xb(-x,-y)$, so that the mapping
is actually well-defined on the quotient $2$-sphere. 

If, in addition,~$\Phi_X$ and~$\Phi_Y$ are positive 
on the closure of~$D_\lambda$, then~$\xb$ is an immersion
away from the half-lattice points.  Moreover, 
the Euler equation~\eqref{eq: Phi Euler} and the smoothness 
of~$\Phi$ near the corner~$(X,Y)=(\lambda_2,\lambda_2)$
implies that~$\Phi_X(\lambda_2,\lambda_2)=\Phi_Y(\lambda_2,\lambda_2)>0$. 
From this, it is not difficult to see that the image~$\xb(T)$
must be a smoothly embedded convex~$2$-sphere. 

In particular, note that holding~$\lambda_2$ fixed and
varying $\lambda_1$ and~$\lambda_3$ gives a $2$-parameter
family of deformations preserving the shape operator 
in a neighborhood of an umbilic point with Hopf index~$\frac12$.

\subsubsection{Examples}\label{sssec: examples II}
In this last section, some interesting examples of these
formulae will be investigated.

\begin{example}[Quadrics]\label{ex: quadrics II}
Consider the \eqref{eq: Phi Euler} solution
\begin{equation}\label{eq: part Phi 1}
\Phi(X,Y) = \frac{2}{\sqrt{-XY}}\,,
\end{equation}
which is defined on the open quadrant~$Y<0<X$.  

Suppose that~$\lambda_1<0<\lambda_3$ and, for simplicity, 
that~$\lambda_2\not=0$. Then the formula~\eqref{eq: x quadrature II} yields
(up to a translation constant)
\begin{equation}\label{eq: part x 1}
\xb(X,Y) =\left[\begin{matrix} 
\ds{\frac1{\lambda_1}\sqrt{\frac{(Y-\lambda_1)(X-\lambda_1)}
                       {(\lambda_2-\lambda_1)(\lambda_3-\lambda_1)(-XY)}}} \\
\ds{\frac1{\lambda_2}\sqrt{\frac{(\lambda_2-Y)(X-\lambda_2)}
                       {(\lambda_2-\lambda_1)(\lambda_3-\lambda_2)(-XY)}}} \\
\ds{\frac1{\lambda_3}\sqrt{\frac{(\lambda_3-Y)(\lambda_3-X)}
                       {(\lambda_3-\lambda_1)(\lambda_3-\lambda_2)(-XY)}}} 
\end{matrix}\,\right]
\end{equation}
defined on the rectangle for which~$\lambda_1<Y<\min\{0,\lambda_2\}$
and~$\max\{0,\lambda_2\}<X<\lambda_3$.
All of these immersions have the shape operator
\begin{equation}\label{eq: part S 1}
S = -\sqrt{-X^3Y}\, \frac{\partial\hfil}{\partial X}\otimes \d X
   +\sqrt{-XY^3}\, \frac{\partial\hfil}{\partial Y}\otimes \d Y\,,
\end{equation}

As is easily checked, the image of~$\xb = (u,v,w)$ 
lies in the hyperboloid of one sheet defined by
\begin{equation}\label{eq: part x image hyperboloid}
\lambda_1\,u^2+\lambda_2\,v^2+\lambda_3\,w^2 
+\frac1{\lambda_1\lambda_2\lambda_3} = 0.
\end{equation}

The case~$\lambda_2=0$ has to be treated separately.  It turns out that
these surfaces are algebraic surfaces of degree~$3$.  Details are left
to the reader.

Finally, consider the \eqref{eq: Phi Euler} solution
\begin{equation}\label{eq: part Phi elliptic 1}
\Phi(X,Y) = \frac{-2}{\sqrt{XY}}\,,
\end{equation}
which is defined on the two open wedges~$0<Y<X$ and~$Y<X<0$. 
Note that~$\Phi_X$ and~$\Phi_Y$ are positive and unequal
on the domain of definition.   Thus, the corresponding mappings~$\xb$
defined on the domains~$D_\lambda$ contained in the domain of~$\Phi$
will be umbilic-free immersions with positive principal curvatures.

If~$0<\lambda_1<\lambda_2<\lambda_3$ or $\lambda_1<\lambda_2<\lambda_3<0$, 
so that~$D_\lambda$ is contained in the domain of~$\Phi$,
then the image of the corresponding immersions~$\xb:D_\lambda\to\E3$
is part of an ellipsoid with three distinct principal axes.
The method of~\S\ref{sssec: torus sphere II} then shows how
this can be used to parametrize the entire ellipsoid.  Note that
the corner~$(\lambda_2,\lambda_2)$ (which lies in the closure of~$D_\lambda$)
gives rise to the four umbilic points that such an ellipsoid possesses.

The cases with~$\lambda_1=0$ give elliptic paraboloids 
while~$\lambda_1<0<\lambda_2<\lambda_3$ gives one sheet 
of an hyperboloid of two sheets.  Details are left to the reader.

\end{example}

\begin{example}[Polynomial solutions]\label{ex: polynomial II}
Consider the quadratic \eqref{eq: Phi Euler} solution
\begin{equation}\label{eq: quadratic part Phi 2}
\Phi(X,Y) = -\frac12(3\,X^2+2\,XY+3\,Y^2)
\end{equation}
defined on the entire half-plane~$Y<X$.
Then the formula~\eqref{eq: x quadrature II} yields
(up to a translation constant)
\begin{equation}\label{eq: part x 2}
\xb(X,Y) =\left[\begin{matrix} 
\ds{\frac{(X+Y+2\lambda_1)\sqrt{(Y-\lambda_1)(X-\lambda_1)}}
         {\sqrt{(\lambda_2-\lambda_1)(\lambda_3-\lambda_1)}}} \\[15pt]
\ds{\frac{(X+Y+2\lambda_2)\sqrt{(\lambda_2-Y)(X-\lambda_2)}}
         {\sqrt{(\lambda_2-\lambda_1)(\lambda_3-\lambda_2)}}} \\[15pt]
\ds{\frac{(X+Y+2\lambda_3)\sqrt{(\lambda_3-Y)(\lambda_3-X)}}
         {\sqrt{(\lambda_3-\lambda_1)(\lambda_3-\lambda_2)}}} \\
\end{matrix}\,\right]
\end{equation}
defined on the rectangle~$D_\lambda$ 
for which~$\lambda_1<Y<\lambda_2<X<\lambda_3$.

By the method of~\S\ref{sssec: torus sphere II}, this can be extended
to a smooth mapping of a $2$-sphere into~$\E3$.
Note that $\xb$ fails to be an immersion along the lines~$Y+3X=0$ and~$X+3Y=0$.
Away from these lines, all of these immersions have the shape operator
\begin{equation}\label{eq: part S 2}
S = \frac{-1}{(3X+Y)}\, \frac{\partial\hfil}{\partial X}\otimes \d X
   +\frac{-1}{(X+3Y)}\, \frac{\partial\hfil}{\partial Y}\otimes \d Y\,,
\end{equation}

It is worth noting that for each positive degree~$d$, 
the equation~\eqref{eq: Phi Euler} has a polynomial solution~$\Phi_d(X,Y)$
that is homogeneous of degree~$d$ and that the solution with this
property is unique up to constant multiples.  The corresponding
surfaces are algebraic, but can be quite complicated.  Nevertheless,
some interesting examples can be found here:

For example, the homogeneous cubic solution
\begin{equation}\label{eq: cubic part Phi 2}
\Phi(X,Y) = -\frac13(5\,X^3+3\,X^2Y+3\,XY^2+5\,Y^3)
\end{equation}
is defined on the entire half-plane~$Y<X$ and its partials~$\Phi_X$
and~$\Phi_Y$ are positive on the entire half-plane.  They are equal
only when~$X+Y=0$ and they vanish on the closed half-plane only at
the point~$(X,Y)=(0,0)$.

As a result, as long as~$\lambda_2\not=0$, the corresponding
immersion~$\xb$ extends to an immersion of the $2$-sphere as a
strictly convex ovaloid.  As long as $\lambda_1+\lambda_2>0$
or~$\lambda_2+\lambda_3<0$, the domain~$D_\lambda$ does not
meet the line~$X+Y=0$, so the only umbilics of the resulting
ovaloid are the four corresponding to~$(\lambda_2,\lambda_2)$.
However, if~$\lambda_1+\lambda_2<0<\lambda_2+\lambda_3$, then 
the domain~$D_\lambda$ meets the line~$X+Y=0$ in a segment, 
which gives rise to a circle of umbilics in the resulting ovaloid.

\end{example}

\begin{example}[Minimal surfaces]\label{ex: minimal surface II}
As a final example, consider the solution
\begin{equation}\label{eq: part Phi 3}
\Phi(X,Y) = 2\,\log(X-Y)
\end{equation}
defined on the entire half-plane~$Y<X$.  Note that, because~$\Phi_X+\Phi_Y=0$,
the resulting surfaces will be minimal surfaces.  In fact, it is 
evident that, up to scaling and the addition of a constant, 
this is the unique solution that satisfies $\Phi_X+\Phi_Y=0$,
and so gives a minimal surface.

The formula~\eqref{eq: x quadrature II} yields
(up to a translation constant)
\begin{equation}\label{eq: part x 3}
\xb(X,Y) =\left[\begin{matrix} 
\ds{\frac{\log\left(\ds\frac{\sqrt{X-\lambda_1}+\sqrt{Y-\lambda_1}}
                            {\sqrt{X-\lambda_1}-\sqrt{Y-\lambda_1}}\right)}
         {\sqrt{(\lambda_2-\lambda_1)(\lambda_3-\lambda_1)}}} \\[15pt]
\ds{-\frac{2\arctan\left(\ds\frac{\sqrt{X-\lambda_2}}
                             {\sqrt{\lambda_2-Y}}\right)}
         {\sqrt{(\lambda_2-\lambda_1)(\lambda_3-\lambda_2)}}} \\[15pt]
\ds{\frac{\log\left(\ds\frac{\sqrt{\lambda_3-Y}-\sqrt{\lambda_3-X}}
                            {\sqrt{\lambda_3-Y}+\sqrt{\lambda_3-X}}\right)}
         {\sqrt{(\lambda_3-\lambda_1)(\lambda_3-\lambda_2)}}}
\end{matrix}\,\right]
\end{equation}
defined on the rectangle~$D_\lambda$.
All of these immersions have the shape operator
\begin{equation}\label{eq: part S 3}
S = \frac{X-Y}{2}\, \frac{\partial\hfil}{\partial X}\otimes \d X
   +\frac{Y-X}{2}\, \frac{\partial\hfil}{\partial Y}\otimes \d Y\,.
\end{equation}

Using~$(u,v,w)$ as coordinates on~$\E3$, the image minimal surface 
satisfies the equation
\begin{equation}\label{eq: implicit minimal II}
\begin{split}
0 &= (\lambda_2-\lambda_3)\,
   \cosh\left(u\,\sqrt{(\lambda_2-\lambda_1)(\lambda_3-\lambda_1)}\,\right)\\
&\qquad -(\lambda_3-\lambda_1)\,
   \cos\left(v\,\sqrt{(\lambda_2-\lambda_1)(\lambda_3-\lambda_2)}\,\right)\\
&\qquad -(\lambda_1-\lambda_2)\,
   \cosh\left(w\,\sqrt{(\lambda_3-\lambda_1)(\lambda_3-\lambda_2)}\,\right).
\end{split}
\end{equation}
Thus, this belongs to the family (investigated
by Weingarten and by Frech\'et, see~\cite[Chapter~II, \S5.2]{Ni}) 
of minimal surfaces that satisfy equations of the form~$f(u)+g(v)+h(w)=0$.
\end{example}

\bibliographystyle{amsplain}

\end{document}